\documentclass[a4paper,10pt]{article}
\usepackage{graphicx,multirow}
\usepackage{amssymb,amsmath,amsthm}
\usepackage{lscape,lipsum,microtype}
\usepackage[margin=2.54cm]{geometry}  
\usepackage{hyperref}

\usepackage{ifluatex}
\ifluatex
 \usepackage{unicode-math}
\fi

\usepackage{latexsym}
\usepackage{amsmath}
\usepackage{amsfonts}
\usepackage{amssymb}
\usepackage{amscd}

\renewcommand{\d}{\delta }
\newcommand{\D }{\Delta }

\renewcommand{\l }{\lambda }

\newcommand{\n }{\nabla }

\newcommand{\intbar}{\mathop{\int\makebox(-13.5,0){\rule[4pt]{.7em}{0.3pt}}%
\kern-6pt}\nolimits}

\newcommand{\be}{\begin{equation}}
\newcommand{\ee}{\end{equation}}
\newcommand{\bes}{\begin{equation*}}
\newcommand{\ees}{\end{equation*}}
\newcommand{\ba}{\begin{eqnarray}}
\newcommand{\ea}{\end{eqnarray}}
\newcommand{\bas}{\begin{eqnarray*}}
\newcommand{\eas}{\end{eqnarray*}}
\newenvironment{pf}{\noindent{\sc Proof}.\enspace}{\rule{2mm}{2mm}\medskip}
\newenvironment{pfn}{\noindent{\sc Proof}}{\rule{2mm}{2mm}\medskip}

\newcommand{\R}{\mathbb{R}}

\newcommand{\Z}{\mathbb{Z}}

\newcommand{\N}{\mathbb{N}}

\author{ Mohammed ALDAWOOD$^a$, Cheikh Birahim NDIAYE$^b$}

\date{}

\title{\bf The Brezis-Nirenberg problem on non-contractible bounded domains of \;$\R^3$}
\begin{document}

\newtheorem{lem}{Lemma}[section]
\newtheorem{pro}[lem]{Proposition}
\newtheorem{thm}[lem]{Theorem}
\newtheorem{rem}[lem]{Remark}
\newtheorem{cor}[lem]{Corollary}
\newtheorem{df}[lem]{Definition}

\maketitle

\begin{center}
{\small

\noindent  $^{a, b}$\; Department of Mathematics Howard University \\  Annex 3, Graduate School of Arts and Sciences, \# 217 \\ DC 20059 Washington, USA.
}

\end{center}


\footnotetext[1]{E-mail addresses: cheikh.ndiaye@howard.edu, mohammed.aldawood@bison.howard.edu\\
\thanks{\\ C. B. Ndiaye was partially supported by NSF grant DMS--2000164.}}

\

\

\begin{center}
{\bf Abstract}
\end{center}
In this paper, we study the Brezis-Nirenberg problem on bounded smooth domains of \;$\R^3$. Using the algebraic topological argument of Bahri-Coron\cite{bc} as implemented in \cite{martndia2}  combined with the Brendle\cite{bre1}-Schoen\cite{sc}'s bubble construction, we solve the problem for non-contractible bounded smooth domains.
 \begin{center}


\bigskip

\centerline{\bf AMS subject classification: 53C21, 35C60, 58J60, 55N10.}

\end{center}


\section{Introduction and statement of the results}
In their seminal paper \cite{brenir}, Brezis and Nirenberg initiated the study of nonlinear elliptic equations of the form
\begin{equation}\label{bnproa}
\left\{
\begin{split}
-\D u +qu&=u^{\frac{n+2}{n-2}}\;\;&\text{in}\;\;\Omega,\\
u&>0\;\;&\text{in}\;\;\Omega,\\
u&=0\;\;&\text{on}\;\;\partial \Omega,
\end{split}
\right.
\end{equation}
where \;$\Omega$\; is a bounded and smooth domain of \;$\R^n$,  $n\geq 3$\; and \;$q$ is a bounded and smooth function defined on \;$\Omega$. In this paper, we revisit the Boundary Value problem (BVP) \eqref{bnproa} in the \;$3$-dimensional case, namely when \;$n=3$. Thus, we will be dealing with the BVP
\begin{equation}\label{eqd3}
\left\{
\begin{split}
-\D u +qu&=u^5\;\;&\text{in}\;\;\Omega,\\
u&>0\;\;&\text{in}\;\;\Omega,\\
u&=0\;\;&\text{on}\;\;\partial \Omega,
\end{split}
\right.
\end{equation}
where \;$\Omega$\; is a smooth bounded domain in \;$\R^3$. It is well known that a necessary condition for the existence of positive solution to \eqref{eqd3} is that the first eigenvalue of $-\D+q$\; under zero Dirichlet boundary condition is positive, see \cite{brezis2}. Morever, we will assume that \;$-\D+q$\; under zero Dirichlet boundary condition verifies the strong maximum principle. Hence the BVP \eqref{eqd3} has a variational structure, since thanks to the strong maximum principle and standard elliptic regularity theory solutions of \eqref{eqd3} can be found by looking at critical points of the Brezis-Nirenberg functional
\begin{equation}\label{eq:dfjq}
J_q(u):=\frac{\left<u, u\right>_q}{(\int_{\Omega}u^6\;dx)^{\frac{1}{3}}}, \;\;\;\;u\in H^{1, +}_0(\Omega):=\{u\in H^1_0(\Omega):\;\;\;\;u\geq 0 \;\;\text{and}\;\;u\neq 0\},
\end{equation}
where 
\begin{equation}\label{scalq}
\left<u, u\right>_q=\int_{\Omega}(|\n u|^2+qu^2)\;dx
\end{equation}
and $H^1_0(\Omega)$\; is the usual Sobolev space of functions which are \;$L^2$-integrable  on \;$\Omega$\; together with their first derivatives and with zero trace on \;$\partial \Omega$.
 \vspace{4pt}

\noindent
Existence of solutions under a Positive Mass type assumption has been obtained in an unpublished work by McLeod as discussed in the work of Brezis\cite{brezis2}. In this work, we use the Barycenter technique of Bahri-Coron\cite{bc} to remove the Positive Mass type assumption of McLeod and replace it by the non-contractibility of the domain. Precisely, we prove the following theorem.
\begin{thm}\label{non-contractible}
Assuming that \;$\Omega\subset \R^3$\;is a non-contractible bounded and smooth domain, $q$\; is a smooth and bounded function defined on \;$\Omega$, the first eigenvalue of the operator \;$-\D+q$\; under zero Dirichlet boundary condition on $\partial \Omega$\; is positive, $-\D+q$\; under zero Dirichlet boundary condition on $\partial \Omega$\;verifies the strong maximum principle and the Green's function \;$G$ of $-\D+q$\; under zero Dirichlet boundary condition on $\partial \Omega$\; defined by \eqref{eqgreen} is positive, then the BVP \eqref{eqd3} has a least one solution. 
\end{thm}
\vspace{4pt}

\noindent
To prove Theorem\ref{non-contractible}, we will use the Algebraic topogical argument of Bahri-Coron\cite{bc} which is possible since as already observed by McLeod (see \cite{brezis2}), the problem under study is a Global one (for the definition of "Gobal" for Yamabe type problems, see \cite{nss}). Indeed, as in \cite{nss}, we will follow the scheme of the Barycenter technique as performed in the work \cite{martndia2} of the second author and Mayer. One of the main difficulty with respect to the works  \cite{martndia2}, and \cite{nss} is the presence of the linear term \;$"qu"$\; and the lack of conformal invariance. Such a difficulty has already been encountered by Bahri-Brezis\cite{bb} on closed Riemannian manifolds. To deal with such a difficulty, Bahri-Brezis\cite{bb} 
have used the bubble construction of Bahri-Coron\cite{bc} recalling that their scheme of the Barycenter technique follows the original one of Bahri-Coron\cite{bc}. However, here we use the Brendle\cite{bre1}-Schoen\cite{sc}' s bubble construction and have to deal with that difficulty in a way different from the work of Bahri-Brezis\cite{bb}.

\section{Notations and preliminaries}
In this section, we fix some notation and discuss some preliminaries. We start with fixing some notation. $\N$\; denotes the  set of non-negative integers and \;$\N^*$\; denotes set of the positive integers. For \;$a\in \R^3$\; and \;$\d>0$, \;$B_a(\d)=B(a, \d)$\; denotes the Euclidean Ball with radius \;$\d$ centered at \;$a$. $\large{1}_{A}$\; denotes the characteristic function of \;$A$. $\n$\; denotes the Euclidean gradient and \;$\D$ denotes the Euclidean Laplacian. All integrations are with respect to \;$dx$\; the standard Lebesque measure on \;$\R^3$\; with\; $x=(x_1, x_2, x_3)$\; the standard coordinate system of \;$\R^3$. $|\cdot|$\; and \;$\left<\cdot, \cdot\right>$\; denote respectively the standard norm and scalar product on \;$\R^3$.  We also use \;$|\cdot|$\; to denote the absolute value on \;$\R$. For \;$E\subset \R^3$\; and \;$p\in \N^*$, $L^p(E)$\; denotes the usual Lebesgue space of order \;$p$\; with respect to \;$dx$. For \;$\chi: \R\longrightarrow \R$\; a smooth function, \;$\chi^{'}$\; and \;$\chi^{''}$ \;denotes respectively  the first derive and second derivative of \;$\chi$. For \;$a\in K\subset \Omega$\;, $K$\; compact, and  \;$0\leq d_1<d_2\leq \infty$, we set \;$\{d_1\leq |x-a|\leq d_2\}=\{x\in \Omega:\;\; d_1\leq |x-a|\leq d_2\}$. To simplify notation, we write\; $d_1\leq |x-a|\leq d_2$\; instead of \;$\{d_1\leq |x-a|\leq d_2\}$\; if there is no possible confusion. Similarly, for \;$0\leq d_1<d_2\leq \infty$, we set \;$\{d_1\leq |y|\leq d_2\}=\{y\in \R^3:\;\; d_1\leq |y|\leq d_2\}$. To simplify notation, we write\; $d_1\leq |y|\leq d_2$\; instead of \;$\{d_1\leq |y|\leq d_2\}$\; if there is no possible confusion.
\vspace{4pt}

\noindent
Next, we introduce the standard bubbles of the variational problem under study. For \;$a\in \R^3$\; and \;$\l>0$, we denote by \;$\d_{a, \l}$\; the standard bubble on  \;$\R^3$, namely
\begin{equation}\label{bubble-R3}
\delta_{a, \;\lambda}(x)=c_0\left(\frac{\lambda}{1+\lambda^2|x-a|^{2}}\right)^{\frac{1}{2}}, \;\;\;\;x\in \R^3,
\end{equation}
where\; $c_0>0$ is such that \;$\d_{a, \l}$\; satisfies
\begin{equation}\label{eqflat}
-\Delta \d_{a, \l}=\d_{a, \l}^{5}\;\;\;\text{on}\;\;\;\R^3.
\end{equation}
We have also the following relations
\begin{equation}\label{s1}
\int_{\R^3}|\n \d_{a, \l}|^2=\int_{\R^3}\d_{a, \l}^6=\int_{\R^3}|\n \d_{0,1}|^2=\int_{\R^3}\d_{0, 1}^6
\end{equation}
and
\begin{equation}\label{s2}
\mathcal{S}=\frac{\int_{\R^3}|\n \d_{a, \l}|^2}{(\int_{\R^3}\d_{a, \l}^6)^{\frac{1}{3}}}, 
\end{equation}
where
\begin{equation}\label{s3}
\mathcal{S}=\inf_{u\in D^1(\R^3), \;u\neq 0}\frac{\int_{\R^3}|\n u|^2}{(\int_{\R^3}u^6)^{\frac{1}{3}}}
\end{equation}
with \;$$D^1(\R^3)=\{u\in L^6(\R^3):\;\;\;\;|\n u| \in L^{2}(\R^3)\}.$$  We set
\begin{equation}\label{c3}
c_3=\int_{\R^3}\left(\frac{1}{1+|y|^2}\right)^{\frac{5}{2}}.
\end{equation}
\vspace{4pt}

\noindent
For \;$a\in \Omega$, let \(G(a,x)\)\; be the unique solution of (see \cite{brezis2})
\begin{equation}\label{eqgreen}
\begin{cases} -\Delta G(a, x)+qG(a, x)=4\pi\delta_a(x),& \;\;\;\;\; x\in\Omega\\\hspace{2.3cm}G(a, x)=0,& \;\;\;\;x\in  \partial\Omega .\end{cases}
\end{equation}
 $G(a, x)$\; satisfies the following estimates
  \begin{equation}\label{estg}
 \left|G(a,x)-\frac{1}{|x-a|}\right|\leq C\;\;\;\text{for} \;\;\;x\neq a\in \Omega,
 \end{equation}
 and
 \begin{equation}\label{estgg}
 \left|\nabla \left(G(a,x)-\frac{1}{|x-a|}\right)\right|\le \frac{C}{|x-a|}\;\; \text{for}\; \;\;x\neq a\in \Omega.
 \end{equation}
 Moreover, under the assumption of Theorem \ref{non-contractible}, we have \;$G>0$\; in \;$\Omega$.
 \vspace{4pt}
 
 \noindent
 Now, let \(\ \chi:\mathbb{R} \to[0,1]\) be a smooth cut-off function satisfying
\begin{equation}\label{chi}
\chi(t)= \left\{ \begin{array}{ll}
         1 & \mbox{if $t \leq  1$}\\
        0 & \mbox{if $t \geq 2$}.\end{array} \right.
        \end{equation}
Using \;$\chi$, for \;$a\in \Omega$, and \;$\delta>0$ and small, we define
\begin{equation}\label{chid}  
        \hspace{-0.5cm}\chi^a_\delta(x)=\chi\left(\frac{|x-a|}{\delta}\right),\;\;\;x\in \Omega.
        \end{equation}       
Moreover, using \;$\chi^a_\delta $\; and the Green's function \;$G(a, \cdot)$, we define the Brendle\cite{bre1}-Schoen\cite{sc}'s bubble 
\begin{equation}\label{uald}
       u_{a, \l, \d}=\chi^a_\delta  \delta_a,_\lambda+(1-\chi^a_\delta)\frac{c_0}{\sqrt{\lambda}}G(a,x).
 \end{equation}
For\;$K\subset \Omega$\; compact, we set \hspace{0.3cm} 
\begin{equation}\label{varrho}
\varrho_0=\varrho_0^K:=\frac{dis(K,\partial\Omega)}{4}>0.  
\end{equation}
Thus, for $\forall a\in K$ and $\forall 0<2\delta<\varrho_0$\; we have 
\begin{equation}\label{ual}
u_{a, \l}:=u_{a, \l, \d}\in H^1_0(\Omega), \;\;\;\text{and}\;\;\; u_{a, \l}>0\;\;\text{in}\;\;\Omega.
\end{equation}
For \;$a_i, a_j\in \Omega$\; and \;$\l_i, \l_j>0$, we define
\begin{equation}\label{varepij}
\varepsilon_{ij}=\left[\frac{1}{\frac{\lambda_i}{\lambda_j}+\frac{\lambda_j}{\lambda_i}+\lambda_i\lambda_jG^{-2}(a_i,a_j)}\right]^{\frac{1}{2}}.
\end{equation}
Moreover, for \;$a_i, a_j\in K$,\; $0<2\delta<\varrho_0$, and \;$\l_i, \l_j>0$, we define
\begin{equation}\label{epij}
\epsilon_{ij}=\int_\Omega u^5_{a_i},_{\lambda_i}u_{a_j},_{\lambda_j}
\end{equation}
and
\begin{equation}\label{eij}
e_{ij}=\int_\Omega(-\Delta+q) u_{a_i},_{\lambda_i}u_{a_j},_{\lambda_j}.
\end{equation}
Using \eqref{eqflat} and \eqref{eqgreen}, we estimate the deficit of \;$u_{a, \l}$\; being a solution of BVP \eqref{eqd3}. 
\begin{lem}\label{c0estimate}
Let \;$K\subset \Omega$\;be  compact and \;$\theta>0$\; be small. Then there exists \;$C>0$\; such that \;$\forall a\in K$,\; $\forall 0<2\delta<\varrho_0$\; and \;$\forall 0<\frac{1}{\l}\leq \theta \d$, we have
\[\left|-\Delta u_{a,\lambda}+q u_{a,\lambda}-u^5_{a,\lambda}\right|\le C\left[\frac{1}{\delta^2\sqrt{\lambda}}\Large{1}_{\{\delta\le|x-a|\le 2\delta\} }+\delta_a,_\lambda\Large{1}_{\{|x-a|\le 2\delta\}}+\delta^5_{a,\lambda} \Large{1}_{\{|x-a|\geq \delta|\}}\right],\]
where \;$\varrho_0$\; is as in \eqref{varrho}.
\end{lem}
\begin{pf}
First of all, to simplify notation, let us set \; $\chi_\delta:=\chi^{a}_{\delta}$,  $G_a(x):=G(a, x)$\; and \;$\bar G_a=c_0G_a$. Then, we have \\\
\[u_a,_\lambda=\chi_\delta \delta_a,_\lambda+(1-\chi_\delta)\frac{\bar G_a}{\sqrt{\lambda}}=\chi_\d\left(\delta_{a, \lambda}-\frac{\bar G_a}{\sqrt{\lambda}}\right)+\frac{\bar G_a}{\sqrt{\lambda}}.\]\\\
This implies
\[\left(-\Delta +q\right)u_a,_\lambda=\left(-\Delta +q\right)\left[\chi_\delta\left(\delta_a,_\lambda-\frac{\bar G_a}{\sqrt{\lambda}}\right)\right]+\frac{\left(-\Delta +q\right)\bar G_a}{\sqrt{\lambda}}.\]\\\
Clearly the lemma is true for \;$x=a$. Now, since for  \;$x\neq a$, we have \;$\left(-\Delta +q\right)\bar G_a=0$,\; then for $x\neq a$\; we get \\
$$\left(-\Delta+q\right)u_a,_\lambda=-\Delta\chi_\delta\left[\delta_a,_\lambda-\frac{\bar G_a}{\sqrt{\lambda}}\right]-2\nabla\chi_\delta \nabla\left[\delta_a,_\lambda-\frac{\bar G_a}{\sqrt{\lambda}}\right]-\chi_\delta\Delta\delta_a,_\lambda+q\chi_\delta\delta_a,_\lambda.$$ This implies (for $x\neq a$)
$$
\left(-\Delta +q\right)u_{a,\lambda}-u_{a, \l}^5=\sum_{i=1}^4 J_i
$$
with 
\begin{equation*}
\begin{split}
&J_1=-\Delta\chi_\delta\left[\delta_a,_\lambda-\frac{\bar G_a}{\sqrt{\lambda}}\right]\\
&J_2=-2\left<\nabla\chi_\delta , \nabla\left[\delta_a,_\lambda-\frac{\bar G_a}{\sqrt{\lambda}}\right]\right>\\
&J_3=q\chi_\delta \delta_a,_\lambda\\
&J_4=-\chi_\delta \Delta \delta_a,_\lambda-u^5_a,_\lambda.
\end{split}
\end{equation*}
Now, we are going to estimate separately each \;$J_i$'s. For \;$J_1$, we first write 
\begin{equation}\label{j1}
J_1=\Delta\chi_\delta\left[\delta_a,_\lambda-\frac{c_0}{\sqrt{\lambda}|x-a|}+\frac{c_0}{\sqrt{\lambda}|x-a|}-\frac{\bar G_a}{\sqrt{\lambda}}\right].
\end{equation} 
Next, using \eqref{bubble-R3} and \eqref{estg}, we derive 
\begin{equation}\label{diff1}
\left|\delta_a,_\lambda-\frac{c_0}{\sqrt{\lambda}|x-a|}\right|\le \frac{C}{\sqrt{\lambda}}
\end{equation}
and 
\begin{equation}\label{diff2}
\left|\frac{c_0}{\sqrt{\lambda}|x-a|}-\frac{\bar G_a}{\sqrt{\lambda}}\right|\le\frac{C}{\sqrt{\lambda}}. 
\end{equation}
For \;$\Delta\chi_\delta$, we have \\\
\begin{equation}\label{gchi}
\nabla\chi_\delta=\chi^{'}\left(\frac{|x-a|}{\delta}\right)\frac{(x-a)}{\delta|x-a|}.
\end{equation}
This implies
\begin{equation}\label{lchi}
\Delta\chi_\delta=\chi^{''}\left(\frac{|x-a|}{\delta}\right)\frac{1}{\delta^2}+2\chi^{'}\left(\frac{|x-a|}{\delta}\right)\frac{1}{\delta|x-a|}.
\end{equation}
Thus, recalling the definition of \;$\chi$ (see \eqref{chi}), we have \eqref{lchi} implies
\begin{equation}\label{flchi}
|\Delta\chi_\delta|\le \frac{C}{\delta^2}\large{1}_{\{\delta\le|x-a|\le 2\delta\}}.
\end{equation}
Hence, combining \eqref{j1}, \eqref{diff1}, \eqref{diff2}, and  \eqref{flchi},  we get 
\begin{equation}\label{estj1}
|J_1|\le \frac{C}{\delta^2\sqrt{\lambda}} \large{1}_{\{\delta\le|x-a|\le 2\delta\}}.
\end{equation}
To estimate $J_2$, we first write 
\begin{equation}\label{j2}
J_2=-2\left<\nabla\chi_\delta , \;\nabla\left[\delta_a,_\lambda-\frac{c_0}{\sqrt{\lambda}|x-a|}+\frac{c_0}{\sqrt{\lambda}|x-a|}-\frac{\bar G_a}{\sqrt{\lambda}}\right]\right>.
\end{equation}
Next,  using \eqref{bubble-R3} and \eqref{estgg}, we derive
\begin{equation}\label{diff3}
\hspace{-2.7cm}\left|\nabla\left[\delta_a,_\lambda-\frac{c_0}{\sqrt{\lambda}|x-a|}\right]\right|\le \frac{C}{\sqrt{\lambda}|x-a|}, 
\end{equation}
and 
\begin{equation}\label{diff4}
\hspace{-4.7cm}\left|\nabla\left[\frac{c_0}{\sqrt{\lambda}|x-a|}-\frac{\bar G_a}{\sqrt{\lambda}}\right]\right|\le \frac{C}{\sqrt{\lambda}|x-a|}.
\end{equation}
On the other hand, using \eqref{gchi} and recalling \eqref{chi}, we derive 
\begin{equation}\label{fgchi}
|\n \chi_\delta|\le \frac{C}{\delta} \large{1}_{\{\delta\le|x-a|\le 2\delta\}}.
\end{equation}
Hence, combining \eqref{j2}-\eqref{fgchi}, we get
\begin{equation}\label{estj2}
|J_2|\le \frac{C}{\delta^2\sqrt{\lambda}} \large{1}_{\{\delta\le|x-a|\le 2\delta\}}.
\end{equation}
For  \;$J_3$, since \;$q$\; is bounded then using \eqref{chi} and \eqref{chid}, we clearly obtain
\begin{equation}\label{fj3}
|J_3|\le C\delta_a,_\lambda \large{1}_{\{|x-a|\le 2\delta\}} .
\end{equation}
Finally to estimate \;$J_4$, we observe that for  \;$|x-a|\le \delta$,  \[\hspace{-5cm}\chi_\delta(x)=1.\]Thus 
\begin{equation}\label{j41}
J_4=-\chi_\d\Delta\delta_a,_\lambda-u^5_a,_\lambda=-\Delta\delta_a,_\lambda-\d^5_a,_\lambda=0
\end{equation}
on $\{|x-a|\leq \d \}$. On  the other hand on \;$\{|x-a|> \delta\}$, we clearly have 
\begin{equation}\label{j42}
|u_a,_\lambda|\le C\delta_a,_\lambda.
\end{equation}
Therefore, \eqref{j41} and \eqref{j42} imply 
\begin{equation}\label{fj4}
\hspace{1.1cm}|J_4|\le\delta^5_a,_\lambda \large{1}_{\{|x-a|\geq \delta\}}.
\end{equation}
Hence, the result follows from  \eqref{estj1}, \eqref{estj2}, \eqref{fj3}, and \eqref{fj4}. 
\end{pf}

\section{PS-sequences and Deformation lemma}
In this section, we recall the analysis of Palais-Smale (PS) sequence for \;$J_q$\; defined by \eqref{eq:dfjq}, see \cite{brezis2}. We also introduce the neighborhood of potential critical points at infinity  of $J_q$ and the associated selection maps. As in other applications of the Barycenter technique of Bahri-Coron\cite{bc}, we also recall the associated Deformation lemma.
\begin{lem}\label{psseq}
Suppose that \;$(u_k)\subset H^{1, +}_0(\Omega)$ is a PS-sequence for \;$J_q$, that is \;$\n J_q(u_k) \rightarrow 0$\; and \;$J_q(u_k)\rightarrow c$ \;up to a subsequence, and $\int_{\Omega} u_k^6=c^{\frac{3}{2}}$, then up to a subsequence, we have have there exists \;$u_{\infty}\geq 0$, an integer \;$p\geq 0$, a sequence of points ${a_{i, k}}\in \Omega, \;\;i=1, \cdots, p$, and a sequence of positive numbers ${\l_{i, k}}, \;\;i=1, \cdots p$,  such that\\
1)\\
$$
-\D u_{\infty}+qu_{\infty}=u_{\infty}^5.
$$
2)\\
$$
||u_k-u_{\infty}-\sum_{i=1}^p u_{a_{i, k}, \l_{i, k}}||_q\longrightarrow 0.
$$
3)\\
$$
J_q(u_k)^{\frac{3}{2}}\longrightarrow J_q(u_{\infty})^{\frac{3}{2}}+pS^{\frac{3}{2}}.
$$
4)\\
For \;$i\neq j=1, \cdots, p$, 
$$ 
\frac{\l_{i, k}}{\l_{j, k}}+\frac{\l_{j, k}}{\l_{i, k}}+\l_{i, k}\l_{j, k}G^{-2}(a_{i, k}, a_{j k})\longrightarrow +\infty
$$
5)\\
For $i=1, \cdots, p$, 
$$
\l_{i, k} dist(a_{i, k}, \partial \Omega)\longrightarrow +\infty,
$$
where \;$||\cdot||_q$\; is the norm associated to the scalar product \;$\left<\cdot, \cdot\right>_q$\; defined by \eqref{scalq}.
\end{lem}
\vspace{6pt}

\noindent
To introduce the neighborhoods of potential critical points at infinity of \;$J_q$, we first fix
\begin{equation}\label{varepsilon0}
\varepsilon_0>0\;\;\; \text{and}\;\;\;\varepsilon_0 \simeq 0.
\end{equation}
 Furthermore, we choose
 \begin{equation}\label{nu0}
 \nu_0>1\;\;\;\text{and}\;\;\;\nu_0\simeq 1.
 \end{equation}
Then for \;$p\in \N^*$,\; and \;$0<\varepsilon\leq \varepsilon_0$, we define \;$V(p, \varepsilon)$\; the \;$(p, \varepsilon)$-neighborhood of potential critical points at infinity of \;$J_q$\; by
\begin{equation*}
\begin{split}
V(p, \varepsilon):=\{u\in H^{1, +}_0(\Omega): &\;\;\exists a_1, \cdots, a_{p}\in \Omega,\;\;\alpha_1, \cdots, \alpha_{p}>0,\;\; \;\l_1, \cdots,\l_{p}>0,  \\&\;\l_i\geq \frac{1}{\varepsilon}\;\;\text{for}\;\;i=1\cdots, p,\;\; \;\;\;\l_i dist (a_i, \partial \Omega)\geq \frac{1}{\varepsilon}\;\;\text{for}\;\;i=1\cdots, p,\;\; \\& \Vert u-\sum_{i=1}^{p}\alpha_i u_{a_i, \l_i}\Vert_q\leq \varepsilon,\;\;\frac{\alpha_i}{\alpha_j}\leq \nu_0\;\;\text{and}\;\;\varepsilon_{i, j}\leq \varepsilon\;\;\text{for}\;\;i\neq j=1, \cdots, p\}.
\end{split}
\end{equation*}
\vspace{6pt}

\noindent
Concerning the sets \;$V(p, \varepsilon)$, for every \;$p\in \N^*$ \;there exists \;$0<\varepsilon_p\leq\varepsilon_0$\; such that for every \;$0<\varepsilon\leq \varepsilon_p$, we have
\begin{equation}\label{eq:mini}
\begin{cases}
\forall u\in V(p, \varepsilon)\;\; \text{the minimization problem}\;\;\min_{B_{\varepsilon}^{p}}\Vert u-\sum_{i=1}^{p}\alpha_iu_{a_i, \l_i}\Vert_q \\
\text{has a solution }\;(\bar \alpha, A, \bar \l)\in B_{\varepsilon}^{p} , \text{which is unique up to permutations,}
\end{cases}
\end{equation}
where \;$B^{p}_{\varepsilon}$\; is defined as 
\begin{equation*}
\begin{split}
B_{\varepsilon}^{p}:=\{(&\bar\alpha=(\alpha_1, \cdots, \alpha_p), A=(a_1, \cdots, a_p), \bar \l=\l_1, \cdots,\l_p))\in \R^{p}_+\times \Omega^p\times (0, +\infty)^{p}\\&\l_i\geq \frac{1}{\varepsilon},\;\l_i dist (a_i, \partial \Omega)\geq \frac{1}{\varepsilon}, i=1, \cdots, p, \;\;\frac{\alpha_i}{\alpha_j}\leq \nu_0\;\; \text{and}\;\ \varepsilon_{i, j}\leq \varepsilon, i\neq j=1, \cdots, p\}.
\end{split}
\end{equation*}
We define the selection map \;$s_p$\; via
\begin{equation*}
s_{p}: V(p, \varepsilon)\longrightarrow (\Omega)^p/\sigma_p\\
:
u\longrightarrow s_{p}(u)=A
\;\,\text{and} \,\;A\;\;\text{is given by}\;\,\eqref{eq:mini}.
\end{equation*}
To state the Deformation Lemma needed for the application of the algebraic topological argument of Bahri-Coron\cite{bc}, we first  set
\begin{equation}\label{dfwp}
W_p:=\{u\in \;:\;J_q(u)\leq (p+1)^{\frac{2}{3}}\mathcal{S}\},
\end{equation}
 for \;$p\in \N$. 
\vspace{4pt}

\noindent
As in \cite{bc}, \cite{martndia2}, and \cite{nss}, we have Lemma \ref{psseq} implies the following Deformation lemma.
\begin{lem}\label{deform}
Assuming that \;$J_q$ \;has no critical points, then for every \;$p\in \N^*$, up to taking \;$\varepsilon_p$\; given by \eqref{eq:mini} smaller, we have that for every\; $0<\varepsilon\leq \varepsilon_p$, the topological pair \;$(W_p,\; W_{p-1})$\; retracts by deformation onto \;$(W_{p-1}\cup A_p, \;W_{p-1})$\; with \;$V(p, \;\tilde \varepsilon)\subset A_p\subset V(p, \;\varepsilon)$\; where \;$0<\tilde \varepsilon<\frac{\varepsilon}{4}$\; is a very small positive real number and depends on \;$\varepsilon$.

\end{lem}
\section{Self-action estimates}
In this section, we derive some sharp estimates needed for application of the Barycenter technique  of Bahri-Coron\cite{bc}. For the numerator of \;$J_q$, we have.
 \begin{lem}\label{num}
 Assuming that \;$K\subset \Omega$\; is compact and \;$\theta>0$\; is small, then there exists \;$C>0$\; such that \;$\forall a\in K$, \;$\forall 0<2\delta<\varrho_0$, and \;$\forall 0<\frac{1}{\l}\leq \theta\delta$, we have 
 \[\int_\Omega (-\Delta+q)u_a,_\lambda u_a,_\lambda \le \int_\Omega u^6_a,_\lambda+\frac{C}{\lambda}\left(1+\delta+\frac{1}{\lambda^2\delta^3}\right).\]

 \end{lem}

\begin{pf}
Setting
 \[I=\int_\Omega (-\Delta+q)u_a,_\lambda u_a,_\lambda,\]
 we get \[I=\int_\Omega u^6_a,_\lambda+\underbrace{\int_{\Omega} \left[(-\Delta+q)u_a,_\lambda -u^5_a,_\lambda\right]u_a,_\lambda}_
{\mbox{$I_1$}}.\] To continue, let us estimate \;$I_1$. Using Lemma \ref{c0estimate}, we get 
\begin{equation*}
\begin{split}
|I_1|&\le \int_\Omega\left|(-\Delta+q)u_a,_\lambda -u^5_a,_\lambda\right|u_a,_\lambda\\&\le \frac{C}{\delta^2\sqrt{\lambda}}\int_\Omega u_{a, \lambda}\Large{1}_{\{\delta\le|x-a|\le 2\delta\}}\\&+C \int_\Omega\delta_a,_\lambda u_a,_\lambda\Large{1}_{\{|x-a|\le 2\delta\}}\\&+C\int_\Omega\delta^5_a,_\lambda u_a,_\lambda\Large{1}_{\{|x-a|\geq \delta\}}.
\end{split}
\end{equation*}
We are going to estimate the three parts of the right hand side  the latter formula. For the first term, we have 
\begin{equation*}
\begin{split}
\int_\Omega u_a,_\lambda\Large{1}_{\{\delta\le|x-a|\le 2\delta\}}&\le C\int_{\delta\le|x-a|\le 2\delta} \left[\frac{\lambda}{1+\lambda^2|x-a|^2}\right]^{\frac{1}{2}}\\&
\le C\int_{\delta\le|x-a|\le 2\delta} \frac{1}{\sqrt{\lambda}|x-a|}\\&\le \frac{C}{\sqrt{\lambda}} \int_{\delta}^{2\delta} r \,dr\\& \le C \frac{\delta^2}{\sqrt{\lambda}}.
\end{split}
\end{equation*}
 For the second term, we obtain
 \begin{equation*}
 \begin{split}
 \int_\Omega\delta_a,_\lambda u_a,_\lambda\Large{1}_{\{|x-a|\le 2\delta\}} &\le C\int_{|x-a|\le 2\delta} \left[\frac{\lambda}{1+\lambda^2|x-a|^2}\right]\\&\le \frac{1}{\lambda} \int_{0}^{2\delta} \,dr\\&\le C\frac{\delta}{\lambda}.
 \end{split}
 \end{equation*}
 Finally for the last term, we get
 \begin{equation*}
 \begin{split}
 \int_\Omega\delta^5_a,_\lambda u_a,_\lambda\Large{1}_{\{|x-a|\geq 2\delta\}} &\le C\int_{|x-a|\geq 2\delta}\delta^6_a,_\lambda\\&\le C\int_{|x-a|\geq 2\delta}\left[\frac{\lambda}{1+\lambda^2|x-a|^2}\right]^{3}\\&\le\frac{C}{\lambda^3}\int_{\{|x-a|\geq 2\delta\}}\frac{1}{|x-a|^6}\\&\le\frac{C}{\lambda^3}\int_{2\delta}^{+\infty} r^{-4}\,dr \\&\le \frac{C}{\lambda^3\delta^3}.
 \end{split}
 \end{equation*}
Thus, collecting all we have
\[|I_1|\le C\left[\frac{1}{\lambda}+\frac{\delta}{\lambda}+\frac{1}{\lambda^3 \delta^3}\right].\]
Hence, we obtain \\
\[\int_\Omega (-\Delta+q)u_a,_\lambda u_a,_\lambda \le \int_\Omega u^6_a,_\lambda+\frac{C}{\lambda}\left(1+\delta+\frac{1}{\lambda^2\delta^3}\right),\]
thereby ending the proof.
\end{pf}
\vspace{4pt}

\noindent
For the denominator of \;$J_q$, we have
\begin{lem}\label{denom}
 Assuming that \;$K\subset \Omega$\; is compact and \;$\theta>0$\; is small, then there exists \;$C>0$\; such that  \;$\forall a\in K$, \;$\forall 0<2\delta<\varrho_0$, and \;$\forall 0<\frac{1}{\l}\leq \theta\delta$, we have
\[\int_\Omega u^6_a,_\lambda=\int_{\R^3} \d^6_a,_\lambda+O\left(\frac{1}{\lambda^3\delta^3}\right).\]
\end{lem}

\begin{pf}
We have
\[\int_\Omega u^6_a,_\lambda=\int_{|x-a|\le \delta} u^6_a,_\lambda+\int_{\delta<|x-a|\le 2\delta} u^6_a,_\lambda+\int_{|x-a|>2\delta} u^6_a,_\lambda.\]
Now, we estimate each term of the right hand side of the latter formula. For the first term, we obtain 
\begin{equation*}
\begin{split}
\int_{|x-a|\le \delta} u^6_a,_\lambda&=\int_{|x-a|\le \delta} \delta^6_a,_\lambda\\&=\int_{\R^3} \delta^6_a,_\lambda-\int_{|x-a|>\delta} \delta^6_a,_\lambda\\&=\int_{\R^3} \delta^6_a,_\lambda+O\left(\frac{1}{\lambda^3\delta^3}\right).
\end{split}
\end{equation*}
For the second term, we derive 
\begin{equation*}
\begin{split}
\int_{\delta<|x-a|\le 2\delta} u^6_a,_\lambda&\le C\int_{\delta\le|x-a|\le 2\delta}\left(\frac{\lambda}{1+\lambda^2|x-a|^2}\right)^{3}\\&\le \frac{C}{\lambda^3}\int_{\delta}^{2\delta} r^{-4}\,dr\\&\le \frac{C}{\lambda^3\delta^3}.
\end{split}
\end{equation*}
For the last term, using \eqref{estg} we get
\begin{equation*}
\begin{split}
\int_{|x-a|\geq 2\delta} u^6_a,_\lambda&=\int_{|x-a|\geq 2\delta} \left(\frac{1}{\sqrt{\lambda}G_a}\right)^{6}\\&=\frac{C}{\lambda^3}\int_{|x-a|\geq2\delta} G^{6}_a\\&
\le\frac{C}{\lambda^3}\int_{|x-a|\geq2\delta} \frac{1}{|x-a|^6}\\&\le\frac{C}{\lambda^3\delta^3}.
\end{split}
\end{equation*}
Therefore, we have \\
\[\int_\Omega u^6_a,_\lambda=\int_{\R^3} \d^6_a,_\lambda+O\left(\frac{1}{\lambda^3\delta^3}\right).\]
\end{pf}
\vspace{4pt}

\noindent
Finally, we derive the \;$J_q$-energy estimate of \;$u_{a, \l}$\; needed for the application of the Barycenter technique of Bahri-Coron\cite{bc}. 
\begin{cor}\label{sharpenergy}
Assuming that \;$K\subset \Omega$\; is compact and \;$\theta>0$\; is small, then there exists \;$C>0$\; such that  \;$\forall a\in K$, \;$\forall 0<2\delta<\varrho_0$, and \;$\forall 0<\frac{1}{\l}\leq \theta\delta$, we have 
$$
J_q(u_{a, \l})\leq \mathcal{S}\left( 1+ C\left[\frac{1}{\l}+\frac{\delta}{\l}+\frac{1}{\delta^3\l^3}\right]\right).
$$
\end{cor}
\begin{pf}
It follows from the properties of \;$\d_{a, \l}$ (see \eqref{s1}-\eqref{s3}), Lemma \ref{num} and Lemma \ref{denom}.
\end{pf}
\section{Interaction estimates}
In this section, we derive sharp inter-action estimates needed for the algebraic topological argument for existence. Recalling \eqref{epij} and \eqref{eij}, we start with the following one relating \;$e_{ij}$\; and \;$\epsilon_{ji}$.
\begin{lem}\label{interact1}
Assuming that \;$K\subset \Omega$\; is compact and \;$\theta>0$\; is small, then there exists \;$C>0$\; such that  \;$\forall a_i, a_j\in K$, \;$\forall 0<2\delta<\varrho_0$, and \;$\forall 0<\frac{1}{\l_i}, \frac{1}{\l_j}\leq \theta\delta$, we have  
\[\int_\Omega u_{a_i},_{\lambda_i}\left|(-\Delta+q)u_{a_j},_{\lambda_j}-u^5_{a_j},_{\lambda_j}\right|\le C\left(\delta+\frac{1}{\lambda^2_j\delta^2}\right)\left(\frac{\lambda_i}{\lambda_j}+\lambda_i\lambda_j|a_i-a_j|^2\right)^{\frac{-1}{2}}.\]
\end{lem}

\begin{pf}
Using Lemma \ref{c0estimate}, we have
\begin{equation}\label{ac0est}
\underbrace{\left|(-\Delta+q)u_{a_j},_{\lambda_j}-u^5_{a_j},_{\lambda_j}\right|}_{\mbox{$L_j$}} \le C\left[\frac{1}{\delta^2\sqrt{\lambda_j}}\Large{1}_{\delta\le|x-a_j|\le 2\delta} +\delta_{a_j},_{\lambda_j}\Large{1}_{|x-a_j|\le 2\delta}+\delta^5_{a_j},_{\lambda_j}\Large{1}_{|x-a_j|\geq \delta}\right].
\end{equation}
On the set \;$\{|x-a_j|\leq 2\d\}$, we have 
$$\delta^2_{a_j, \lambda_j}(x)\geq c_0^2\left(\frac{\lambda_j}{1+4\lambda^2_j \delta^2}\right)\geq\frac{c^2_0}{\lambda_j \delta^2}\left[1+O(\frac{1}{\lambda^2_j \delta^2})\right]\geq\frac{1}{2}\frac{c_0^2}{\lambda_j\delta^2}.$$
This implies $$\frac{1}{\sqrt{\lambda_j} \delta}\leq \frac{\sqrt{2}}{c_0^2}\delta_{a_j, \lambda_j}.$$
Thus, we get 
\[L_j\le C\left(1+\frac{1}{\delta}\right)\delta_{a_j},_{\lambda_j} \Large{1}_{|x-a|\le 4\delta}+C\delta^5_{a_j},_{\lambda_j} \Large{1}_{|x-a|\geq \frac{\delta}{2}}, \]
where \;$L_j$\; is as in \eqref{ac0est}.
Hence, we obtain 
\begin{equation}\label{intest}
\begin{split}
\int_\Omega u_{a_i},_{\lambda_i} L_j \le &C\left(1+\frac{1}{\delta}\right) \underbrace{\int_{|x-a_j|\le4\delta}\left(\frac{\lambda_j}{1+\lambda_j^2|x-a_j|^2}\right)^{\frac{1}{2}}\left(\frac{\lambda_i}{1+\lambda_i^2|x-a_i|^2}\right)^{\frac{1}{2}}}_{\mbox{$I_1$}}
\\&+C\underbrace{\int_{|x-a_j|\geq \frac{\delta}{2}}\left(\frac{\lambda_j}{1+\lambda_j^2|x-a_j|^2}\right)^{\frac{5}{2}}\left(\frac{\lambda_i}{1+\lambda_i^2|x-a_i|^2}\right)^{\frac{1}{2}}}_{\mbox{$I_2$}}.
\end{split}
\end{equation}
Now, we estimate $I_1$ as follows.
\[\hspace{-3.6cm}I_1=\int_{|x-a_j|\le4\delta}\left(\frac{\lambda_j}{1+\lambda_j^2|x-a_j|^2}\right)^{\frac{1}{2}}\left(\frac{\lambda_i}{1+\lambda_i^2|x-a_i|^2}\right)^{\frac{1}{2}}\]\[=\underbrace{\int_{(2|x-a_i|\le\frac{1}{\lambda_j}+|a_i-a_j|)\cap( |x-a_j|\le 4\delta)}\left(\frac{\lambda_j}{1+\lambda_j^2|x-a_j|^2}\right)^{\frac{1}{2}}\left(\frac{\lambda_i}{1+\lambda_i^2|x-a_i|^2}\right)^{\frac{1}{2}}}_{\mbox{$I^1_1$}}\]
\[+\underbrace{\int_{(2|x-a_i|>\frac{1}{\lambda_j}+|a_i-a_j|)\cap( |x-a_j|\le 4\delta)}\left(\frac{\lambda_j}{1+\lambda_j^2|x-a_j|^2}\right)^{\frac{1}{2}}\left(\frac{\lambda_i}{1+\lambda_i^2|x-a_i|^2}\right)^{\frac{1}{2}}}_{\mbox{$I^2_1$}}.\]
To continue, we first estimate $I^1_1$. Indeed, using triangle inequality we have 
\[I^1_1\le C\int_{|x-a_i|\le 8\delta}\left(\frac{\lambda_j}{1+\lambda_j^2|a_i-a_j|^2}\right)^{\frac{1}{2}}\left(\frac{\lambda_i}{1+\lambda_i^2|x-a_i|^2}\right)^{\frac{1}{2}}\]
\[\hspace{-1.5cm}\le C\frac{\sqrt{\frac{\lambda_j}{\lambda_i}}}{\left(1+\lambda^2_j|a_i-a_j|^2\right)^{\frac{1}{2}}}\int_{|x-a_i|\le 8\delta} \frac{1}{|x-a_i|}.\]
This implies
\begin{equation}\label{i11}
I^1_1=O\left(\delta^2\left(\frac{\lambda_i}{\lambda_j}+\lambda_i\lambda_j|a_i-a_j|^2\right)^{\frac{-1}{2}}\right).
\end{equation}
For $I^2_1$, we derive
\[I^2_1\le C \int_{|x-a_j|\le 4\delta}\left(\frac{\lambda_j}{1+\lambda_j^2|x-a_j|^2}\right)^{\frac{1}{2}}\left(\frac{\lambda_i}{1+\lambda_i^2|a_i-a_i|^2}\right)^{\frac{1}{2}}\left(\frac{\lambda_j}{\lambda_i}\right)\]
\[\hspace{-2.5cm}\le C\frac{\sqrt{\frac{\lambda_i}{\lambda_j}}\left(\frac{\lambda_j}{\lambda_i}\right)}{\left(1+\lambda^2_j|a_i-a_j|^2\right)^{\frac{1}{2}}}\int_{|x-a_j|\le 4\delta} \frac{1}{|x-a_j|}.\]
Thus for $I^2_1$, we obtain
\begin{equation}\label{i12}
I^2_1=O\left(\delta^2\left(\frac{\lambda_i}{\lambda_j}+\lambda_i\lambda_j|a_i-a_j|^2\right)^{\frac{-1}{2}}\right).
\end{equation}
 Hence, combining \eqref{i11} and \eqref{i12}, we get 
 \begin{equation}\label{i1}
I_1=O\left(\delta^2\left(\frac{\lambda_i}{\lambda_j}+\lambda_i\lambda_j|a_i-a_j|^2\right)^{\frac{-1}{2}}\right).
\end{equation}
Next, let us estimate $I_2$. For this, we first write
\[\hspace{-3.6cm}I_2=\int_{|x-a_j|\geq \frac{\delta}{2}}\left(\frac{\lambda_j}{1+\lambda_j^2|x-a_j|^2}\right)^{\frac{5}{2}}\left(\frac{\lambda_i}{1+\lambda_i^2|x-a_i|^2}\right)^{\frac{1}{2}}\]
\[=\underbrace{\int_{\{2|x-a_i|\le\frac{1}{\lambda_j}+|a_i-a_j|\}\cap(|x-a_j|\geq\frac{\delta}{2})}\left(\frac{\lambda_j}{1+\lambda_j^2|x-a_j|^2}\right)^{\frac{5}{2}}\left(\frac{\lambda_i}{1+\lambda_i^2|x-a_i|^2}\right)^{\frac{1}{2}}}_{\mbox{$I^1_2$}}\]
\[+\underbrace{\int_{(2|x-a_i|>\frac{1}{\lambda_j}+|a_i-a_j|)\cap(|x-a_j|\geq\frac{\delta}{2})}\left(\frac{\lambda_j}{1+\lambda_j^2|x-a_j|^2}\right)^{\frac{5}{2}}\left(\frac{\lambda_i}{1+\lambda_i^2|x-a_i|^2}\right)^{\frac{1}{2}}}_{\mbox{$I^2_2$}}.\]
 Setting \;$\mathcal{D}=\{2|x-a_i|\le\frac{1}{\lambda_j}+|a_i-a_j|\}\cap\{|x-a_j|\geq \frac{\d}{2}\}$, we estimate \;$I^1_2$\; as follows
\begin{equation*}
\begin{split}
I^1_2 &\le \frac{C}{ \lambda_j^{\frac{-5}{2}}}\int_{\mathcal{D}}\left(\frac{1}{1+\lambda_j^2|a_i-a_j|^2}\right)^{\frac{3}{2}}\left(\frac{\lambda_i}{1+\lambda_i^2|x-a_i|^2}\right)^{\frac{1}{2}}\left(\frac{1}{1+\lambda_j^2|a_i-a_j|^2}\right)\\
&\le C\left(\frac{1}{\sqrt{\lambda_i}}\right)\left(\frac{1}{\lambda_j^2 \delta^2}\right) \frac{\lambda_j^{\frac{5}{2}}}{\left(1+\lambda_j^2|a_i-a_j|^2\right)^{\frac{3}{2}}}\int_{2|x-a_i|\le\frac{1}{\lambda_j}+|a_i-a_j|}\frac{1}{|x-a_i|}\\
&\le C \left(\sqrt{\frac{\lambda_j}{\lambda_i}}\right)\left(\frac{1}{\delta^2}\right) \frac{1}{\left(1+\lambda_j^2|a_i-a_j|^2\right)^{\frac{3}{2}}}\left(1+\lambda_j^2|a_i-a_j|^2\right)^2\left(\frac{1}{\lambda_j^2}\right)\\
&\le C  \left(\sqrt{\frac{\lambda_j}{\lambda_i}}\right)\left(\frac{1}{\lambda_j^2 \delta^2}\right) \frac{1}{\left(1+\lambda_j^2|a_i-a_j|^2\right)^{\frac{1}{2}}}.
\end{split}
\end{equation*}
This implies 
\begin{equation}\label{i21}
I^2_1=O\left(\frac{1}{\lambda^2_j\delta^2}\left(\frac{\lambda_i}{\lambda_j}+\lambda_i\lambda_j|a_i-a_j|^2\right)^{\frac{-1}{2}}\right).
\end{equation}
Next, we estimate \;$I^2_2$\; as follows
\begin{equation*}
\begin{split}
I^2_2&\le C\frac{\lambda_j}{\sqrt{\lambda_i}}\int_{|x-a_j|\geq\frac{\delta}{2}}\left(\frac{1}{1+\lambda_j^2|a_i-a_j|^2}\right)^{\frac{1}{2}}\left(\frac{\lambda_j}{1+\lambda_j^2|x-a_j|^2}\right)^{\frac{5}{2}}\\
&\le C\frac{\lambda_j}{\sqrt{\lambda_i}}\left(\frac{1}{1+\lambda_j^2|a_i-a_j|^2}\right)^{\frac{1}{2}}\frac{1}{\lambda_j^{\frac{5}{2}}}\int_{|x-a_j|\geq\frac{\delta}{2}}\frac{1}{|x-a_j|^5}\\
&\le C\left(\sqrt{\frac{\lambda_j}{\lambda_i}}\right)\left(\frac{1}{\lambda_j^2 \delta^2}\right) \frac{1}{\left(1+\lambda_j^2|a_i-a_j|^2\right)^{\frac{1}{2}}}.
\end{split}
\end{equation*}
This  gives 
\begin{equation}\label{i22}
I^2_2=O\left(\frac{1}{\lambda^2_j\delta^2}\left(\frac{\lambda_i}{\lambda_j}+\lambda_i\lambda_j|a_i-a_j|^2\right)^{\frac{-1}{2}}\right).
\end{equation}
Therefore, using \eqref{i21} and \eqref{i22}, we obtain 
 \begin{equation}\label{i2}
I_2=O\left(\frac{1}{\lambda^2_j\delta^2}\left(\frac{\lambda_i}{\lambda_j}+\lambda_i\lambda_j|a_i-a_j|^2\right)^{\frac{-1}{2}}\right).
\end{equation}
Hence, combining \eqref{intest},  \eqref{i1} and \eqref{i2}, we get
\[\int_\Omega u_{a_i},_{\lambda_i} L_j \le C\left(\delta+\frac{1}{\lambda^2_j\delta^2}\right)\left(\frac{\lambda_i}{\lambda_j}+\lambda_i\lambda_j|a_i-a_j|^2\right)^{\frac{-1}{2}},\]
thereby ending the proof of the lemma.
\end{pf}
\vspace{4pt}

\noindent
Clearly Lemma \ref{interact1} implies the following sharp interaction-estimate relating \;$e_{ij}$, $\epsilon_{ij}$, and \;$\varepsilon_{ij}$ (for their definitions, see \eqref{varepij}-\eqref{eij}).
\begin{cor}\label{interact2}
Assuming that \;$K\subset \Omega$\; is compact, \;$\theta>0$ is small,  and $\mu_0>0$ is small, then \;$\forall a_i, a_j\in K$, \;$\forall 0<2\delta<\varrho_0$, and \;$\forall 0<\frac{1}{\l_j}\leq\frac{1}{\l_i}\leq \theta\delta$\; such that \;$\varepsilon_{ij}\leq \mu_0$, we have
$$
e_{ij}=\epsilon_{ij}+O\left(\d+\frac{1}{\l_i^2\d^2}\right)\varepsilon_{ij}.
$$
\end{cor}
\vspace{4pt}

\noindent
The next lemma provides a sharp inter-action estimate relating \;$\epsilon_{ji}$\; and \;$\varepsilon_{ij}$.
\begin{lem}\label{interact3}
Assuming that \;$K\subset \Omega$\; is compact,  \;$\theta>0$ is small, and \;$\mu_0>0$\; is small, then  \;$\forall a_i, a_j\in K$, \;$\forall 0<2\delta<\varrho_0$, and \;$\forall 0<\frac{1}{\l_j}\leq\frac{1}{\l_i}\leq \theta\delta$\; such that \;$\varepsilon_{ij}\leq \mu_0$,  we have 
\[\epsilon_{ji}=c_0^6c_3\varepsilon_{ij}\left[\left(1+O(\delta+\frac{1}{\lambda^2_i\delta^2})\right)\left(1+o_{\varepsilon_{ij}}(1)+O(\varepsilon_{ij}^2(\d^{-2}+\log\varepsilon_{ij}^{-1}))\right)+O\left(\varepsilon_{ij}^2\frac{1}{\delta^6}\right)\right],\]
where  \;$c_0$ is as in \eqref{bubble-R3} and \;$c_3$\; is as in \eqref{c3} .
\end{lem}
\begin{pf}
By definition, we have
\[u_{a_i},_{\lambda_i}=\chi_\delta\delta_{a_i},_{\lambda_i}+(1-\chi_\delta)\frac{c_0}{\sqrt{\lambda}}G_{a_i}, \]  with $\chi_\d:=\chi_\d^{a_i}$. On the other hand, by definition of the standard bubble $\d_{a, \l}$, we have 
\[\chi_\delta \delta_{a_i},_{\lambda_i}=c_0\chi_\delta\left[\frac{\lambda_i}{1+\lambda^2_i G^{-2}_{a_i}\frac{|x-a|^2}{G^{-2}_{a_i}}}\right]^{\frac{1}{2}}.\]
Next, for $|x-a_i|\le 2\delta$, we have \\
\begin{equation*}
\begin{split}
1+\lambda^2_i G^{-2}_{a_i}\frac{|x-a_i|^2}{G^{-2}_{a_i}}&=1+\lambda^2_i G^{-2}_{a_i}\left(1+O(\delta)\right)\\&=1+\lambda^2_i G^{-2}_{a_i}+O\left(\lambda^2_i G^{-2}_{a_i} \delta\right)\\
&=\left(1+\lambda^2_i G^{-2}_{a_i}\right)\left[1+O\left(\frac{\lambda^2_i G^{-2}_{a_i} \delta}{1+\lambda^2_i G^{-2}_{a_i}}\right)\right]\\&=\left(1+\lambda^2_i G^{-2}_{a_i}\right)\left[1+O(\delta)\right].
\end{split}
\end{equation*}
So, for $\chi_\delta \delta_{a_i},_{\lambda_i}$ we have
\begin{equation}\label{part1}
\chi_\delta \delta_{a_i},_{\lambda_i}=c_0\chi_\d\left[\frac{\lambda_i}{\left(1+\lambda^2_i G^{-2}_{a_i}\right)\left[1+O(\delta)\right]}\right]^{\frac{1}{2}}
=c_0 \chi_\delta\left[1+O(\delta)\right]\left[\frac{\lambda_i}{1+\lambda^2_i G^{-2}_{a_i}}\right]^{\frac{1}{2}}.
\end{equation}
We have also 
\[c_0(1-\chi_\delta)\left[\frac{\lambda_i}{1+\lambda^2_i G^{-2}_{a_i}}\right]^{\frac{1}{2}}=(1-\chi_\delta)\frac{c_0}{\sqrt{\lambda_i}}G_{a_i}\left[\frac{1}{1+\lambda^{-2}_i G^{2}_{a_i}}\right]^{\frac{1}{2}}.\] Since on $\{|x-a_i|\geq \d\}$, we have
\[ \frac{1}{1+\lambda^{-2}_iG^{2}_{a_i}}=1+O\left(\frac{G^2_{a_i}}{\lambda^2_i}\right)=1+O\left(\frac{1}{\lambda^2_i \delta^2}\right),\] then
we get
\[c_0(1-\chi_\delta)\left[\frac{\lambda_i}{1+\lambda^2_i G^{-2}_{a_i}}\right]^{\frac{1}{2}}=(1-\chi_\delta)\frac{c_0}{\sqrt{\lambda_i}}G_{a_i}\left(1+O\left(\frac{1}{\lambda^2_i \delta^2}\right)\right).\]
This implies 
\begin{equation}\label{part2}
(1-\chi_\delta)\frac{c_0}{\sqrt{\lambda}}G_{a_i}=c_0(1-\chi_\delta)\left[\frac{\lambda}{1+\lambda^2_i G^{-2}_{a_i}}\right]^{\frac{1}{2}}\left(1+O\left(\frac{1}{\lambda^2_i \delta^2}\right)\right).
\end{equation}
Thus, combining \eqref{part1} and \eqref{part2}, we get
\[u_{a_i},_{\lambda_i}=c_0\left[\left(1+O(\delta)\right)\chi_\delta+(1-\chi_\delta)\left(1+O\left(\frac{1}{\lambda^2_i \delta^2}\right)\right)\right]\left[\frac{\lambda}{1+\lambda^2_i G^{-2}_{a_i}}\right]^{\frac{1}{2}}.\]
Hence, we obtain
\begin{equation}\label{partf}
u_{a_i, \l_i}=c_0\left[1+O(\delta)+O\left(\frac{1}{\lambda^2_i\delta^2}\right)\right]\left[\frac{\lambda}{1+\lambda^2_i G^{-2}_{a_i}}\right]^{\frac{1}{2}}.
\end{equation}
Now, we  are going to use \eqref{partf} to achieve our goal. First of all, we write\\
\[\int_\Omega u^5_{a_j},_{\lambda_j}u_{a_i},_{\lambda_i}=\int_{B(a_j, \delta)} u^5_{a_j},_{\lambda_j}u_{a_i},_{\lambda_i}+\int_{\Omega-B(a_j, \delta)} u^5_{a_j},_{\lambda_j}u_{a_i},_{\lambda_i}.\]
For the second term in the right hand side of  the latter formula, we have
\begin{equation*}
\begin{split}
\int_{\Omega-B(a_j,\delta)} u^5_{a_j},_{\lambda_j}u_{a_i},_{\lambda_i}&\le C\int_{\Omega-B(a_j,\delta)}\left(\frac{1}{\lambda_j}\right)^{\frac{5}{2}}\left(\frac{1}{\delta}\right)^5u_{a_i},_{\lambda_i}\\&\le C \left(\frac{1}{\lambda_j}\right)^{\frac{5}{2}}\left(\frac{1}{\delta}\right)^5\int_{\Omega-B(a_j,\delta)}u_{a_i},_{\lambda_i}\\&
\le C\left(\frac{1}{\lambda_j\d^2}\right)^{\frac{5}{2}}\left[\int_{\Omega-(B(a_j,\delta)\cup B(a_i,\delta))}u_{a_i},_{\lambda_i}+\int_{(\Omega-B(a_j,\delta))\cap B(a_i,\delta)}u_{a_i},_{\lambda_i}\right]\\&
\le C\left(\frac{1}{\lambda_j}\right)^{\frac{5}{2}}\left(\frac{1}{\delta}\right)^6\frac{1}{\sqrt{\lambda_i}}+C\left(\frac{1}{\lambda_j}\right)^{\frac{5}{2}}\left(\frac{1}{\delta}\right)^5\int_{B(a_i, \delta)}\left(\frac{\lambda_i}{1+\lambda^2_i|x-a_i|^2}\right)^{\frac{1}{2}}\\&
\le C\left(\frac{1}{\lambda_j}\right)^{\frac{5}{2}}\left(\frac{1}{\delta}\right)^6\frac{1}{\sqrt{\lambda_i}}+C\delta^2\left(\frac{1}{\lambda_j}\right)^{\frac{5}{2}}\left(\frac{1}{\delta}\right)^5\frac{1}{\sqrt{\lambda_i}}\\&\le C\left(\frac{1}{\lambda_j}\right)^{\frac{5}{2}}\left(\frac{1}{\delta}\right)^6\frac{1}{\sqrt{\lambda_i}}\left(1+\delta^3\right).
\end{split}
\end{equation*}
Thus, we  get
\begin{equation}\label{inout}
\int_\Omega u^5_{a_j},_{\lambda_j}u_{a_i},_{\lambda_i}=\int_{B(a_j,\delta)} u^5_{a_j},_{\lambda_j}u_{a_i},_{\lambda_i}+O\left(\frac{1}{\lambda^{\frac{5}{2}}_j\sqrt{\lambda_i}\delta^6}\right).
\end{equation}
For the first term in the right hand side of \eqref{inout} the latter formula, using \eqref{partf} we have
\begin{equation*}
\begin{split}
\int_{B(a_j, \delta)} u^5_{a_j},_{\lambda_j}u_{a_i},_{\lambda_i}&=c^6_0\int_{B(a_j,\delta)}\left(\frac{\lambda_j}{1+\lambda^2_j|x-a_j|^2}\right)^{\frac{5}{2}}\left[1+O(\delta)+O\left(\frac{1}{\lambda^2_i\delta^2}\right)\right]\left[\frac{\lambda_i}{1+\lambda^2_i G^{-2}_{a_i}}\right]^{\frac{1}{2}}\\
&=c^6_0\left[1+O(\delta)+O\left(\frac{1}{\lambda^2_i\delta^2}\right)\right]\int_{B(a_j,\delta)}\left(\frac{\lambda_j}{1+\lambda^2_j|x-a_j|^2}\right)^{\frac{5}{2}}\left[\frac{\lambda_i}{1+\lambda^2_i G^{-2}_{a_i}}\right]^{\frac{1}{2}}\\&
=c^6_0 \frac{1}{\sqrt{\lambda_j}}\left[1+O\left(\d+\frac{1}{\lambda^2_i\delta^2}\right)\right]\int_{B(0,\lambda_j\delta)}\left(\frac{1}{1+|y|^2}\right)^{\frac{5}{2}}\left[\frac{\lambda_i}{1+\lambda^2_i G^{-2}_{a_i}\left(\frac{y}{\lambda_j}+a_j\right)}\right]^{\frac{1}{2}}\\&=c^6_0 \left[1+O\left(\d+\frac{1}{\lambda^2_i\delta^2}\right)\right]\int_{B(0,\lambda_j\delta)}\left(\frac{1}{1+|y|^2}\right)^{\frac{5}{2}}\left[\frac{\lambda_i}{\frac{\l_j}{\l_i}+\lambda_i \l_jG^{-2}_{a_i}\left(\frac{y}{\lambda_j}+a_j\right)}\right]^{\frac{1}{2}}.
\end{split}
\end{equation*}
Recalling that $\lambda_i\le\lambda_j$, then  for \ $\varepsilon_{ij}\sim 0 $, we have\\\
\noindent
1) Either $\varepsilon^{-2}_{ij}\sim \lambda_i\lambda_jG^{-2}_{a_i}(a_j),$\\

\noindent
2) or $\varepsilon^{-2}_{ij}\sim\frac{\lambda_j}{\lambda_i}$.\\
To continue, let
\[\mathcal{A}=\left\{\left|\frac{y}{\lambda_j}\right|\le \epsilon \:G^{-1}_{a_i}(a_j)\right\}\cap B(0,\delta\lambda_j)\cup \left\{\left|\frac{y}{\lambda_j}\right|\le\frac{\epsilon}{\lambda_i}\right\}\cap B(0,\delta\lambda_j),\]
with \;$\epsilon>0$\; very small.  Then  by Taylor expansion on \;$\mathcal{A}$, we have
\begin{equation*}
\begin{split}
\left[\frac{\lambda_j}{\lambda_i}+\lambda_i\lambda_jG^{-2}_{a_i}\left(\frac{y}{\lambda_j}+a_j\right)\right]^{\frac{-1}{2}}&=\left[\frac{\lambda_j}{\lambda_i}+\lambda_i\lambda_jG^{-2}_{a_i}(a_j)\right]^{\frac{-1}{2}}
\\&+\left[\left(\frac{-1}{2}\nabla G^{-2}_{a_i}(a_j)\lambda_i y\right)\right]\left[\frac{\lambda_j}{\lambda_i}+\lambda_i\lambda_jG^{-2}_{a_i}\left(\frac{y}{\lambda_j}+a_j\right)\right]^{\frac{-3}{2}}\\&+O\left[\left(\frac{\lambda_i}{\lambda_j}\right)|y|^2\right]\left[\frac{\lambda_j}{\lambda_i}+\lambda_i\lambda_jG^{-2}_{a_i}\left(\frac{y}{\lambda_j}+a_j\right)\right]^{\frac{-3}{2}}.
\end{split}
\end{equation*}
Thus, we have 
\begin{equation}\label{in}
\int_{B(a_j,\delta)} u^5_{a_j},_{\lambda_j}u_{a_i},_{\lambda_i}=c_0^6\left[1+O(\delta)+O\left(\frac{1}{\lambda^2_i\delta^2}\right)\right]\left(\sum_{i=1}^{4} I_i\right),
\end{equation}
with
\[\hspace{-2.7cm}I_1=\left[\frac{\lambda_j}{\lambda_i}+\lambda_i\lambda_jG^{-2}_{a_i}(a_j)\right]^{\frac{-1}{2}}\int_{\mathcal{A}}\left(\frac{1}{1+|y|^2}\right)^{\frac{5}{2}},\]
\[I_2=\left[\frac{\lambda_j}{\lambda_i}+\lambda_i\lambda_jG^{-2}_{a_i}(a_j)\right]^{\frac{-3}{2}}\int_{\mathcal{A}}\left(\frac{1}{1+|y|^2}\right)^{\frac{5}{2}}\left[\nabla G^{-2}_{a_i}(a_j)\lambda_i y\right],\]
\[I_3=\left[\frac{\lambda_j}{\lambda_i}+\lambda_i\lambda_jG^{-2}_{a_i}(a_j)\right]^{\frac{-3}{2}}\int_{\mathcal{A}}\left(\frac{1}{1+|y|^2}\right)^{\frac{5}{2}}O\left[\left(\frac{\lambda_i}{\lambda_j}\right)|y|^2\right],\]
and
\[I_4=\int_{B(0,\lambda_j\delta)-\mathcal{A}}\left(\frac{1}{1+|y|^2}\right)^{\frac{5}{2}}\left[\frac{\lambda_j}{\lambda_i}+\lambda_i\lambda_jG^{-2}_{a_i}\left(\frac{y}{\lambda_j}+a_j\right)\right]^{\frac{-1}{2}}.\]
Now, let us estimate estimate \;$I_1$. We have 
\[\hspace{-1cm}I_1=\left[\frac{\lambda_j}{\lambda_i}+\lambda_i\lambda_jG^{-2}_{a_i}(a_j)\right]^{\frac{-1}{2}}\left[c_3+\int_{\R^3-\mathcal{A}}\left(\frac{1}{1+|y|^2}\right)^{\frac{5}{2}}\right],\]
where  \;$c_3$ is as in \eqref{c3}. On the other hand, we have
\[ \int_{\R^3-\mathcal{A}}\left(\frac{1}{1+|y|^2}\right)^{\frac{5}{2}} \le \int_{\R^3-B(0,\delta\lambda_j)}\left(\frac{1}{1+|y|^2}\right)^{\frac{5}{2}}+\int_{\R^3-B(0,\lambda_j\epsilon G^{-1}_{a_i}(a_i,a_j))}\left(\frac{1}{1+|y|^2}\right)^{\frac{5}{2}}\]\\
if $\varepsilon^{-2}_{ij}\sim \lambda_i\lambda_jG^{-2}_{a_i}(a_j)$, and 
\[ \int_{\R^3-\mathcal{A}}\left(\frac{1}{1+|y|^2}\right)^{\frac{5}{2}} \le \int_{\R^3-B(0,\delta\lambda_j)}\left(\frac{1}{1+|y|^2}\right)^{\frac{5}{2}}+\int_{\Re^3-B(0,\epsilon\frac{\lambda_j}{\lambda_i})}\left(\frac{1}{1+|y|^2}\right)^{\frac{5}{2}}\]\\
if $\varepsilon^{-2}_{ij}\sim\frac{\lambda_j}{\lambda_i}.$ We have 
\\
\[\hspace{-4cm}\int_{\R^3-B(0,\delta\lambda_j)}\left(\frac{1}{1+|y|^2}\right)^{\frac{5}{2}}=O\left(\frac{1}{\lambda_j^2\delta^2}\right).\] Moreover, if $\varepsilon^{-2}_{ij}\sim \lambda_i\lambda_jG^{-2}_{a_i}(a_j)$, then
\[\hspace{-1.5cm}\int_{\R^3-B(0,\lambda_j\epsilon G^{-1}_{a_i}(a_j))}\left(\frac{1}{1+|y|^2}\right)^{\frac{5}{2}}=O\left(\frac{1}{\lambda^2_j\epsilon^2 G^{-2}_{a_i}(a_j)}\right)\]\[\hspace{5.5cm}=O\left(\frac{1}{\lambda_j\lambda_i G^{-2}_{a_i}(a_j)}\right)=O\left(\varepsilon^2_{ij}\right).\]\\
Furthermore if \;$\varepsilon^{-2}_{ij}\sim\frac{\lambda_j}{\lambda_i}$, then\\
\[\hspace{-6cm}\int_{\R^3-B(0,\epsilon\frac{\lambda_j}{\lambda_i})}\left(\frac{1}{1+|y|^2}\right)^{\frac{5}{2}}=O\left(\varepsilon^2_{ij}\right).\]
This implies \\
\[\int_{\R^3-\mathcal{A}}\left(\frac{1}{1+|y|^2}\right)^{\frac{5}{2}}=O\left(\varepsilon^2_{ij}+\frac{1}{\l_j^2\delta^2}\right)=O\left(\varepsilon^2_{ij}+\varepsilon^2_{ij}\frac{1}{\delta^2}\right)=O\left(\varepsilon^2_{ij}\frac{1}{\delta^2}\right).\]
Thus, we get\ 
\[I_1=\left[\frac{\lambda_j}{\lambda_i}+\lambda_i\lambda_jG^{-2}_{a_i}(a_j)\right]^{\frac{-1}{2}}\left[c_3+O\left(\varepsilon^2_{ij}\frac{1}{\delta^2}\right)\right]\]\[\hspace{-1cm}=\varepsilon_{ij}\left(1+o_{\varepsilon_{ij}}(1)\right)\left[c_3+O\left(\varepsilon^2_{ij}\frac{1}{\delta^2}\right)\right].\]\[\hspace{1cm}\]
Hence, we obtain
\begin{equation}\label{i1fin}
I_1=c_3\varepsilon_{ij}\left[1+o_{\epsilon_{ij}}(1)+O\left(\varepsilon^2_{ij}\frac{1}{\delta^2}\right)\right].
\end{equation}
By symmetry, we have
 \begin{equation}\label{i2fin}
 I_2=0.
 \end{equation}
Next, to estimate \;$I_3$\; we first observe that
\begin{equation*}
\begin{split}
\int_{\mathcal{A}}\frac{|y|^2}{\left(1+|y^2|\right)^{\frac{5}{2}}}&\le \int_{B(0,\epsilon\lambda_jG^{-1}_{a_i}(a_j))-B(0, 1)}\frac{|y|^2}{\left(1+|y^2|\right)^{\frac{5}{2}}}+\int_{B(0,\epsilon \frac{\lambda_j}{\lambda_i})-B(0, 1)}\frac{|y|^2}{\left(1+|y^2|\right)^{\frac{5}{2}}}+O(1)\\&
=O\left(\log(\epsilon\lambda_jG^{-1}_{a_i}(a_j))+\log(\epsilon \frac{\lambda_j}{\lambda_i})\right)+O(1).
\end{split}
\end{equation*}
Thus, we have 
\begin{equation*}
\begin{split}
I_3&=\varepsilon^3_{ij}\left(\frac{\lambda_i}{\lambda_j}\right)\left(1+o_{\varepsilon_{ij}}(1)\right)\left[O\left(\log(\epsilon\lambda_jG^{-1}_{a_i}(a_j))+\log(\epsilon \frac{\lambda_j}{\lambda_i})\right)+O(1)\right]\\&
=\varepsilon^3_{ij}\left(1+o_{\varepsilon_{ij}}(1)\right)\left[O\left(\log(\l_i\lambda_jG^{-2}_{a_i}(a_j))+\log(\frac{\lambda_j}{\lambda_i})\right)+O(1)\right].
\end{split}
\end{equation*}
Hence, we obtain
\begin{equation}\label{i3fin}
I_3=O\left(\varepsilon^3_{ij}\log(\varepsilon_{ij}^{-1})\right).
\end{equation}
Finally, we estimate $I_4$ as follows. \\
If $\varepsilon^{-2}_{ij}\sim\frac{\lambda_j}{\lambda_i}$, then\\
\begin{equation}\label{part1i4}
I_4\le C\varepsilon_{ij}\int_{B(0,\lambda_j\delta)-\mathcal{A}}\left(\frac{1}{1+|y|^2}\right)^{\frac{5}{2}}\le C\varepsilon_{ij}\left(\frac{\lambda_j}{\lambda_i}\right)^{-2}\le C \varepsilon^5_{ij}.
\end{equation}
If $\varepsilon^{-2}_{ij}\sim \lambda_i\lambda_jG^{-2}_{a_i}(a_j)$, then we argue as follows. In case  \;$\left|a_i-a_j\right|\geq2\delta$, since 
\[G_{a_i}\left(\frac{y}{\lambda_j}+a_j\right)\le C\delta^{-1}\]
 for \;$y\in B(0, \l_j\d)$,  then we have 
\begin{equation*}
\begin{split}
I_4&\le C\int_{B(0,\lambda_j\delta)-\mathcal{A}}\left(\frac{1}{1+|y|^2}\right)^{\frac{5}{2}} \frac{1}{\sqrt{\lambda_i\lambda_j}\delta}\\&\le\frac{C}{\sqrt{\lambda_i\lambda_j}\delta}\left(\frac{1}{\lambda^2_jG_{a_i}^{-2}(a_j)}\right)\\&\le \frac{C}{\lambda_i\lambda_jG^{-2}_{a_i}(a_j)}\frac{G^{-1}_{a_i}(a_j)}{\sqrt{\lambda_i\lambda_j}G_{a_i}^{-1}(a_j)\d}\\&\le C\varepsilon^2_{ij} \varepsilon_{ij}\frac{1}{\delta} .
\end{split}
\end{equation*}
Thus, when\; $|a_i-a_j|\geq 2\d$\; we have 
\begin{equation}\label{part2i4}
I_4=O\left(\varepsilon^3_{ij}\frac{1}{\delta}\right).
\end{equation}
In case \; $\left|a_i-a_j\right|<2\delta$, we first observe that 
$$
B(0, \l_j\d)\setminus \mathcal{A}\subset A_1\cup A_2
$$
with 
$$
A_1=\left\{\epsilon\lambda_jG^{-1}_{a_i}(a_j)\le|y|\le E\lambda_jG^{-1}_{a_i}(a_j)\right\}
$$
and
$$
A_2=\left\{E\lambda_jG^{-1}_{a_i}(a_j)\le|y|\le \lambda_j\delta\right\},
$$
where \;$0<\epsilon<E$.
 Thus, we have 
\begin{equation}\label{i4ine}
I_4\le I^1_4+I^2_4,
\end{equation}
with
\begin{equation*}
I_4^1=\int_{A_1} \left(\frac{1}{1+|y|^2}\right)^{\frac{5}{2}}\left[\frac{\lambda_j}{\lambda_i}+\lambda_i\lambda_jG^{-2}_{a_i}\left(\frac{y}{\lambda_j}+a_j\right)\right]^{\frac{-1}{2}}
\end{equation*}
and
\begin{equation*}
I_4^2=\int_{A_2} \left(\frac{1}{1+|y|^2}\right)^{\frac{5}{2}}\left[\frac{\lambda_j}{\lambda_i}+\lambda_i\lambda_jG^{-2}_{a_i}\left(\frac{y}{\lambda_j}+a_j\right)\right]^{\frac{-1}{2}}.
\end{equation*}
We estimate \;$I_4^1$\; as follows:
\begin{equation*}
\begin{split}
\vspace{0.5cm}I_4^1&\le C\left[1+\lambda_j^2G^{-2}_{a_i}(a_j)\right]^{\frac{-5}{2}}\int_{|y|\le E\lambda_jG^{-1}_{a_i}(a_j)}\left[\frac{\lambda_j}{\lambda_i}+\lambda_i\lambda_jG^{-2}_{a_i}\left(\frac{y}{\lambda_j}+a_j\right)\right]^{\frac{-1}{2}}\\&
\le C\left[1+\lambda_j^2G^{-2}_{a_i}(a_j)\right]^{\frac{-5}{2}}\left(\frac{\lambda_i}{\lambda_j}\right)^{\frac{1}{2}}\int_{|y|\le E\lambda_jG^{-1}_{a_i}(a_j)}\left[1+\lambda^2_iG^{-2}_{a_i}\left(\frac{y}{\lambda_j}+a_j\right)\right]^{\frac{-1}{2}}\\&
\le C\left[1+\lambda_j^2G^{-2}_{a_i}(a_j)\right]^{\frac{-5}{2}}\left(\frac{\lambda_i}{\lambda_j}\right)^{\frac{1}{2}}\int_{|y|\le E\lambda_jG^{-1}_{a_i}(a_j)}\left[1+\lambda^2_i\left|\frac{y}{\lambda_j}+a_i-a_j\right|^2\right]^{\frac{-1}{2}}\\&
\le C\left[1+\lambda_j^2G^{-2}_{a_i}(a_j)\right]^{\frac{-5}{2}}\left(\frac{\lambda_i}{\lambda_j}\right)^{\frac{1}{2}}\left(\frac{\lambda_j}{\lambda_i}\right)^{3}\int_{|z|\le \bar E\lambda_iG^{-1}_{a_i}(a_j)}\left[\frac{1}{1+|z|^2}\right]^{\frac{1}{2}}\\&
\le C \left[\frac{\lambda_i}{\lambda_j}+\lambda_i\lambda_jG^{-2}_{a_i}(a_j)\right]^{\frac{-5}{2}}\left(\lambda^2_iG^{-2}_{a_i}(a_j)\right)
\\&
\le C \varepsilon^5_{ij}\left(\lambda_i\lambda_j G^{-2}_{a_i}(a_j)\right),
\end{split}
\end{equation*}
where $\bar E$ is a positive constant. So we obtain
\begin{equation}\label{i14fin}
I^1_4=O\left(\varepsilon^3_{ij}\right).
\end{equation}
For $I_4^2$, we have
\begin{equation*}
\begin{split}
I^2_4&=\int_{B_2}\left(\frac{1}{1+|y|^2}\right)^{\frac{5}{2}}\left[\frac{\lambda_j}{\lambda_i}+\lambda_i\lambda_jG^{-2}_{a_i}\left(\frac{y}{\lambda_j}+a_j\right)\right]^{\frac{-1}{2}}\\&
\le C\int_{|y|\geq E\lambda_jG^{-1}_{a_i}(a_j)}\left(\frac{1}{1+|y|^2}\right)^{\frac{5}{2}}\left[\frac{\lambda_j}{\lambda_i}+\lambda_i\lambda_jG^{-2}_{a_i} (a_j )\right]^{\frac{-1}{2}}
\\&\le C\left[\frac{\lambda_j}{\lambda_i}+\lambda_i\lambda_jG^{-2}_{a_i} (a_j )\right]^{\frac{-1}{2}}\left(\frac{1}{\lambda^2_jG^{-2}_{a_i} (a_j )}\right).
\end{split}
\end{equation*}
This implies
 \begin{equation}\label{i24fin}
 I^2_4=O\left(\varepsilon^3_{ij}\right).
 \end{equation}
Thus, combining \eqref{i4ine}-\eqref{i24fin}, we have that if \;$|a_i-a_j|<2\d$, then 
\begin{equation}\label{i4pcclose}
I_4=O\left(\varepsilon^3_{ij}\right).
\end{equation}
Now, using \eqref{part2i4} and \eqref{i4pcclose}, we infer that  in case \;$\varepsilon_{i, j}^{-2}\simeq \l_i\l_jG_{a_i}^{-2}(a_j)$, 
\begin{equation}\label{i4pc}
I_4=O\left(\varepsilon_{i, j}^3\frac{1}{\d}\right).
\end{equation}
Finally combining \eqref{part1i4} and \eqref{i4pc}, we get
\[I_4=O\left(\varepsilon^3_{ij}\frac{1}{\delta}\right).\]
Collecting all  we get 
\begin{equation}\label{in}
\int_{B(a_j,\lambda_j\delta)} u^5_{a_j},_{\lambda_j}u_{a_i},_{\lambda_i}=c_0^6\left[1+O\left(\d+\frac{1}{\lambda^2_i\delta^2}\right)\right]\left[c_3\varepsilon_{ij}\left(1+o_{\varepsilon_{i, j}}(1)+O(\varepsilon_{i, j}^2(\d^{-2}+\log\varepsilon_{ij}^{-1}))\right)\right].
\end{equation}
Therefore using \eqref{inout} and \eqref{in}, we arrive to 
\begin{equation*}
\begin{split}
\int_\Omega u^5_{a_j},_{\lambda_j}u_{a_i},_{\lambda_i}=&c_0^6\left[1+O\left(\d+\frac{1}{\lambda^2_i\delta^2}\right)\right]\left[c_3\varepsilon_{ij}\left(1+o_{\varepsilon_{i, j}}(1)+O(\varepsilon_{i, j}^2(\d^{-2}+\log\varepsilon_{ij}^{-1}))\right)\right]\\&+O\left(\frac{1}{\lambda^{\frac{5}{2}}_j\sqrt{\lambda_i}\delta^6}\right).
\end{split}
\end{equation*}
Thus, we have
\begin{equation*}
\begin{split}
\int_\Omega u^5_{a_j},_{\lambda_j}u_{a_i},_{\lambda_i}=&c_0^6\left[1+O\left(\d+\frac{1}{\lambda^2_i\delta^2}\right)\right]\left[c_3\varepsilon_{ij}\left(1+o_{\varepsilon_{i, j}}(1)+O(\varepsilon_{i, j}^2(\d^{-2}+\log\varepsilon_{ij}^{-1}))\right)\right]\\&+O\left(\varepsilon_{ij}^3\frac{1}{\delta^6}\right).
\end{split}
\end{equation*}
Therefore, we obtain
\begin{equation}\label{eqfinal}
\begin{split}
\int_\Omega u^5_{a_j},_{\lambda_j}u_{a_i},_{\lambda_i}=&c_0^6c_3\varepsilon_{i,j}\left[\left(1+O\left(\d+\frac{1}{\lambda^2_i\delta^2}\right)\right)\left(1+o_{\varepsilon_{i, j}}(1)+O(\varepsilon_{i, j}^2(\d^{-2}+\log\varepsilon_{ij}^{-1}))\right)\right]\\&+O\left(\varepsilon_{ij}^3\frac{1}{\delta^6}\right).
\end{split}
\end{equation}
Hence the result follows from \eqref{epij} and \eqref{eqfinal}.
\end{pf}
\vspace{4pt}

\noindent
Clearly switching the index \;$i$\; and \;$j$\; in Lemma \ref{interact3}, we have the following corollary which is equivalent to Lemma \ref{interact3}. We decide to present the following corollary, because its form suits more our presentation of the Barycenter technique of Bahri-Coron\cite{bc} which follows the works \cite{martndia2} and \cite{nss}.
\begin{cor}\label{interact4}
Assuming that \;$K\subset \Omega$\; is compact,\;$\theta>0$ is small, and \;$\mu_0$\; is small, then \;$\forall a_i, a_j\in K$, \;$\forall 0<2\delta<\varrho_0$, and \;$\forall 0<\frac{1}{\l_i}\leq\frac{1}{\l_j}\leq \theta\delta$\; such that \;$\varepsilon_{ij}\leq \mu_0$,  we have 
\[\epsilon_{ij}=c_0^6c_3\varepsilon_{ij}\left[\left(1+O\left(\d+\frac{1}{\lambda^2_j\delta^2}\right)\right)\left(1+o_{\varepsilon_{ij}}(1)+O(\varepsilon_{ij}^2(\d^{-2}+\log\varepsilon_{ij}^{-1})\right)+O\left(\varepsilon_{ij}^2\frac{1}{\delta^6}\right)\right].\]
\end{cor}
\vspace{4pt}

\noindent
We present  now some sharp high-order inter-action estimates needed for the application of the algebraic topological argument for existence. We start with the following balanced high-order inter-action estimate.
\begin{lem}\label{interact5}
Assuming that \;$K\subset \Omega$\; is compact, \;$\theta>0$ is small, and $\mu_0$ is small, then \;$\forall a_i, a_j\in K$, \;$\forall 0<2\delta<\varrho_0$, and \;$\forall 0<\frac{1}{\l_j},\frac{1}{\l_i}\leq \theta\delta$\; such that \;$\varepsilon_{ij}\leq \mu_0$, we have
\[\int_\Omega u^{3}_{a_i},_{\lambda_i} u^{3}_{a_j},_{\lambda_j}=O\left(\varepsilon^{3}_{ij}\delta^{-6}\log\left(\frac{1}{\varepsilon_{ij}\d}\right)\right).\]
\end{lem}
\begin{pf}
By symmetry, we can assume without loss of generality (w.l.o.g) that  \;$\lambda_j\le\lambda_i$. Thus we have 
 \\
1) Either $\varepsilon^{-2}_{ij}\sim \lambda_i\lambda_jG^{-2}_{a_i}(a_j)$\\
2) Or $\varepsilon^{-2}_{ij}\sim\frac{\lambda_i}{\lambda_j}$.\\\\
Now, if $|a_i-a_j|\geq2\delta$, then we have 
\begin{equation}\label{i1}
\begin{split}
I&=\int_\Omega u^{3}_{a_i},_{\lambda_i} u^{3}_{a_j},_{\lambda_j}\\&\le C\int_{B(a_i,\delta)}\left(\frac{\lambda_i}{1+\lambda^2_i|x-a_i|^2}\right)^{\frac{3}{2}} \left(\frac{\lambda_j}{1+\lambda^2_jG^{-2}_{a_j}(x)}\right)^{\frac{3}{2}}\\&
+\frac{C}{\lambda_i^{\frac{3}{2}}\delta^3}\int_{B(a_j,\delta)}\left(\frac{\lambda_j}{1+\lambda^2_j|x-a_j|^2}\right)^{\frac{3}{2}}+\frac{C}{\delta^6}\left(\frac{1}{\lambda_i\lambda_j}\right)^{\frac{3}{2}}\\&\le\underbrace{C\int_{B(0,\lambda_i\delta)}\left(\frac{1}{1+|y|^2}\right)^{\frac{3}{2}}\left(\frac{1}{\frac{\lambda_i}{\lambda_j}+\lambda_i\lambda_jG^{-2}_{a_j}\left(a_i+\frac{y}{\lambda_i}\right)}\right)^{\frac{3}{2}}}_{\mbox{$I_1$}}\\&
+\frac{C}{\delta^3}\left(\frac{1}{\lambda_i\lambda_j}\right)^{\frac{3}{2}}\int_{B(0,\lambda_j\delta)}\left(\frac{1}{1+|y|^2}\right)^{\frac{3}{2}}+\frac{C}{\delta^6}\left(\frac{1}{\lambda_i\lambda_j}\right)^{\frac{3}{2}}\\&
\le C I_1+\frac{C}{\delta^3}\left(\frac{1}{\lambda_i\lambda_j}\right)^{\frac{3}{2}}\left[\log(\lambda_j\delta)+C\right]+\frac{C}{\delta^6}\left(\frac{1}{\lambda_i\lambda_j}\right)^{\frac{3}{2}}\\&
\le C I_1+\frac{C}{\delta^6}\left(\frac{1}{\lambda_i\lambda_j}\right)^{\frac{3}{2}}\log(\lambda_j).
\end{split}
\end{equation}
Now, we estimate $I_1$ as follows. \\
If  \;$\varepsilon^{-2}_{ij}\sim\frac{\lambda_i}{\lambda_j}$ , then we get
\[I_1\le C \varepsilon^3_{ij}\left[\log(\lambda_i\delta)+C\right].\]
So, for \;$I$\; we have
\begin{equation*}
\begin{split}
I&\le C \varepsilon^3_{ij}\left[\log(\lambda_i\delta)+C\right]+\frac{C}{\delta^6}\varepsilon^3_{ij}\log(\lambda_i\lambda_j)
\\&\le \frac{C}{\delta^6}\varepsilon^3_{ij}\log(\varepsilon^{-2}_{ij}G^2_{a_i}(a_j))\\&
=O\left(\frac{\varepsilon^3_{ij}\log(\varepsilon^{-1}_{ij}\delta^{-1})}{\delta^6}\right).
\end{split}
\end{equation*}
If \;$\varepsilon^{-2}_{ij}\sim \lambda_i\lambda_jG^{-2}_{a_i}(a_j)$, then we get
\begin{equation*}
I_1\le \frac{C}{\delta^3}\left(\frac{1}{\lambda_i\lambda_j}\right)^{\frac{3}{2}}\left[\log(\lambda_i\delta)+C\right].
\end{equation*}
So, for \;$I$\; we have
\begin{equation*}
I\le \frac{C}{\delta^6}\left(\frac{1}{\lambda_i\lambda_j}\right)^{\frac{3}{2}}\log(\lambda_i\lambda_j).
\end{equation*}
This implies 
\begin{equation*}
I\le \frac{C}{\delta^6}\varepsilon^{3}_{ij}\log(\varepsilon^{-2}_{ij}G^2_{a_i}(a_j)).
\end{equation*} 
Hence,  for \;$|a_i-a_j|\geq 2\d$, we obtain
\begin{equation}\label{estfar}
I=O\left(\frac{\varepsilon^3_{ij}\log(\varepsilon^{-1}_{ij}\delta^{-1})}{\delta^6}\right).
\end{equation}
On the other hand, arguing as above, if $|a_i-a_j|<2\delta$\; then we have also 
\begin{equation*}
\begin{split}
I&\le I_1+\frac{C}{\delta^3}\left(\frac{1}{\lambda_i\lambda_j}\right)^{\frac{3}{2}}\log(\lambda_j)+\frac{C}{\delta^6}\left(\frac{1}{\lambda_i\lambda_j}\right)^{\frac{3}{2}}\\&
\le I_1+\frac{C}{\delta^6}\left(\frac{1}{\lambda_i\lambda_j}\right)^{\frac{3}{2}}\log(\lambda_j\lambda_i),
\end{split}
\end{equation*}
where \;$I_1$\; is as in \eqref{i1}. Thus, if  $\varepsilon_{i,j}^{-2}\simeq \frac{\l_i}{\l_j}$\; then
$$
I\le I_1+\frac{C}{\delta^6}\left(\frac{\l_j}{\lambda_i}\right)^{\frac{3}{2}}\frac{1}{\l_j^3}\left[\log(\frac{\lambda_i}{\l_j})+\log( \l_j^2)\right].
$$
This implies
$$
I\leq  I_1+\frac{C}{\d^6}\varepsilon_{ij}^3\log (\varepsilon_{ij}^{-1}).
$$
Next, if \;$\varepsilon_{i,j}^{-2}\simeq \l_i\l_j G^{-2}_{a_i}(a_j)$\; then we get
\begin{equation*}
\begin{split}
 I&\le I_1+\frac{C}{\delta^6}\left(\frac{1}{\lambda_i\lambda_jG^{-2}_{a_i}(a_j)}\right)^{\frac{3}{2}}\left[\log(\l_i\l_j G^{-2}_{a_i}(a_j)+\log(G^{2}_{a_i}(a_j)))\right]G^{-3}_{a_i}(a_j)\\&
\le I_1+\frac{C}{\delta^6}\varepsilon^3_{ij}\log(\varepsilon^{-1}_{ij}).
\end{split}
\end{equation*}
Now, to continue, we are going to estimate $I_1$. For this, we start by defining the following sets:
\[A_1=\left\{|y|\le \epsilon\lambda_i \sqrt{G^{-2}_{a_j}(a_i)+\frac{1}{\lambda^2_j}}\right\}\]
\[\hspace{4.5cm}A_2=\left\{\epsilon\lambda_i \sqrt{G^{-2}_{a_j}(a_i)+\frac{1}{\lambda^2_j}}\le |y|\le E \lambda_i\sqrt{G^{-2}_{a_j}(a_i)+\frac{1}{\lambda^2_j}}\right\}\]
\[\hspace{1.5cm}A_3=\left\{E\lambda_i \sqrt{G^{-2}_{a_j}(a_i)+\frac{1}{\lambda^2_j}}\le |y|\le 4\lambda_i\delta\right\},\]
with \;$0<\epsilon<E<\infty$. Clearly by  the definition of \;$I_1$ (see \eqref{i1}), we have
\[I_1\le \int_{A_1}L_{ij}+\int_{A_2}L_{ij}+\int_{A_3}L_{ij},\]
where 
\[L_{ij}=\left(\frac{1}{1+|y|^2}\right)^{\frac{3}{2}}\left(\frac{1}{\frac{\lambda_i}{\lambda_j}+\lambda_i\lambda_jG^{-2}_{a_j}\left(a_i+\frac{y}{\lambda_i}\right)}\right)^{\frac{3}{2}}.\]  For \;$\int_{A_1}L_{ij} $, we have 
\begin{equation*}\begin{split}
\int_{A_1}L_{ij} &\le C \varepsilon^{3}_{ij}\int_{A_1}\left(\frac{1}{1+|y|^2}\right)^{\frac{3}{2}}\\
&\le C\varepsilon^{3}_{ij} \log\left(\sqrt{\frac{\lambda_i}{\lambda_j}}\sqrt{\lambda_i\lambda_jG^{-2}_{a_j}(a_i)+\frac{\lambda_i}{\lambda_j}}\right)\\&
\le C\varepsilon^{3}_{ij} \log(\varepsilon^{-1}_{ij}).
\end{split}
\end{equation*}
For \;$\int_{A_2}L_{i,j}$, we have
\begin{equation*}
\begin{split}
\int_{A_2}L_{ij}&\le C \left(\frac{1}{\left(\frac{\lambda_i}{\lambda_j}\right)^2+\lambda^2_iG^{-2}_{a_j}(a_i)}\right)^{\frac{3}{2}}\int_{A_2}\left(\frac{1}{\frac{\lambda_i}{\lambda_j}+\lambda_i\lambda_jG^{-2}_{a_j}\left(a_i+\frac{y}{\lambda_i}\right)}\right)^{\frac{3}{2}}\\&
\le C \left(\frac{\lambda_j}{\lambda_i}\right)^{\frac{3}{2}}\varepsilon^{3}_{ij}\int_{|y|\leq E\l_i\sqrt{G_{a_j}^{-2}(a_i)+\frac{1}{\l_j^2}} }\left(\frac{1}{\frac{\lambda_i}{\lambda_j}+\frac{\lambda_j}{\lambda_i}\left|y+\lambda_i(a_i-a_j)\right|^2}\right)^{\frac{3}{2}}\\&
\le C \left(\frac{\lambda_j}{\lambda_i}\right)^{\frac{3}{2}}\varepsilon^{3}_{ij}\int_{|y|\leq \bar E\l_i\sqrt{G_{a_j}^{-2}(a_i)+\frac{1}{\l_j^2}} }\left(\frac{1}{\frac{\lambda_i}{\lambda_j}+\frac{\lambda_j}{\lambda_i}\left|y\right|^2}\right)^{\frac{3}{2}}\\&
\le C \left(\frac{\lambda_j}{\lambda_i}\right)^{3}\left(\frac{\lambda_j}{\lambda_i}\right)^{-3}\varepsilon^{3}_{ij}\int_{|y|\leq \bar E\l_j\sqrt{G_{a_j}^{-2}(a_i)+\frac{1}{\l_j^2}}}\left(\frac{1}{1+|y|^2}\right)^{\frac{3}{2}}\\&
\le C\varepsilon^{3}_{ij} \log(\varepsilon^{-1}_{ij}).
\end{split}
\end{equation*}
For \;$\int_{A_3}L_{i,j}$, we have 
\begin{equation*}
\begin{split}
\int_{A_3}L_{ij}&\le\int_{A_3}\left(\frac{1}{1+|y|^2}\right)^{\frac{3}{2}}\left(\frac{1}{\frac{\lambda_i}{\lambda_j}+\frac{\lambda_j}{\lambda_i}|y|^2}\right)^{\frac{3}{2}}\\&
\le C \left(\frac{\lambda_i}{\lambda_j}\right)^{\frac{3}{2}}\int_{A_3}\frac{1}{|y|^6}\\&
\le C \left(\frac{\lambda_i}{\lambda_j}\right)^{\frac{3}{2}} \left(\frac{1}{\lambda_i^2G^{-2}_{a_j}(a_i)+\left(\frac{\lambda_i}{\lambda_j}\right)^2}\right)^{\frac{3}{2}}\\&
\le C \varepsilon^{3}_{ij}.
\end{split}
\end{equation*}
Therefore, we have  
$$
I_1\leq C\varepsilon^{3}_{ij} \log(\varepsilon^{-1}_{ij}).
$$
This implies for \;$|a_i-a_j|<2\d$, we have 
$$
I=O\left(\frac{\varepsilon^{3}_{ij}}{\delta^6}\log(\varepsilon^{-1}_{ij})\right).
$$
Hence, combining with the estimate for $|a_i-a_j|\geq 2\d$ (see \eqref{estfar}), we have 
\[\int_\Omega u^{3}_{a_i},_{\lambda_i} u^{3}_{a_j},_{\lambda_j}= O\left(\frac{\varepsilon^{3}_{ij}}{\delta^6}\log(\varepsilon^{-1}_{ij}\d^{-1})\right).\]
\end{pf}
\vspace{4pt}

\noindent
Finally, we present  a sharp unbalanced high-order inter-action estimate needed for the application of the Barycenter technique of Bahri-Coron\cite{bc}.
\begin{lem}\label{interact6}
Assuming that \;$K\subset \Omega$\; is compact, \;$\theta>0$\; is small, and \;$\mu_0$\; is small, then \;$\forall a_i, a_j\in K$, \;$\forall 0<2\delta<\varrho_0$, and \;$\forall 0<\frac{1}{\l_i}\leq \frac{1}{\l_j}\leq \theta\delta$\; such that \;$\varepsilon_{ij}\leq \mu_0$, we have
\[\int_\Omega u^{\alpha}_{a_i},_{\lambda_i} u^{\beta}_{a_j},_{\lambda_j}=O\left(\frac{\varepsilon_{ij}^{\beta}}{\d^6}\right).\]
where \(\alpha+\beta=6\)  and \(\alpha>3>\beta>1.\)
\end{lem}
\begin{pf}
Let \;$\hat{\alpha}=\frac{1}{2}\alpha$ and $\hat{\beta}=\frac{1}{2}\beta$. Then  we have $\hat{\alpha}+\hat{\beta}=3$. Now, since $\lambda_j\le \lambda_i$, then we have\\
1) Either \;$\varepsilon^{-2}_{ij}\sim \lambda_i\lambda_jG^{-2}_{a_i}(a_j)$\\\\
2) Or \;$\varepsilon^{-2}_{ij}\sim\frac{\lambda_i}{\lambda_j}$.\\\
To continue, we write
\[\int_\Omega u^{\alpha}_{a_i, \lambda_i} u^{\beta}_{a_j, \lambda_j}=\underbrace{\int_{B_{a_i}(\delta)} u^{\alpha}_{a_i, \lambda_i} u^{\beta}_{a_j, \lambda_j}}_{\mbox{$I_1$}}+\underbrace{\int_{\Omega-B_{a_i}(\delta)} u^{\alpha}_{a_i, \lambda_i} u^{\beta}_{a_j, \lambda_j}}_{\mbox{$I_2$}}\]
and estimate $I_1$ and $I_2$. For $I_2$, we have 
\begin{equation*}\begin{split}
I_2&= \int_{\left(\Omega-B_{a_i}(\delta)\right)\cap B_{a_j}(\delta)} u^{\alpha}_{a_i},_{\lambda_i} u^{\beta}_{a_j},_{\lambda_j}+\int_{\Omega-\left(B_{a_i}(\delta)\cup B_{a_j}(\delta)\right)} u^{\alpha}_{a_i},_{\lambda_i} u^{\beta}_{a_j},_{\lambda_j}\\&
\le C\int_{\left(\Omega-B_{a_i}(\delta)\right)\cap B_{a_j}(\delta)} \left(\frac{\lambda_i}{1+\lambda^2_iG^{-2}_{a_i}(x)}\right)^{\hat{\alpha}}\left(\frac{\lambda_j}{1+\lambda^2_j|x-a_j|^2}\right)^{\hat{\beta}}\\&
+C\int_{\Omega-\left(B_{a_i}(\delta)\cup B_{a_j}(\delta)\right)}\left(\frac{\lambda_i}{1+\lambda^2_iG^{-2}_{a_i}(x)}\right)^{\hat{\alpha}}\left(\frac{\lambda_j}{1+\lambda^2_jG^{-2}_{a_j}(x)}\right)^{\hat{\beta}}\\&
\le \frac{C}{\lambda^{\hat{\alpha}}_i\lambda^{3-\hat{\beta}}_j\delta^{\alpha}}\int_{B_{0}(\l_j\delta)}\left(\frac{1}{1+|y|^2}\right)^{\hat{\beta}}+\frac{C}{\lambda^{\hat{\alpha}}_i\lambda^{\hat{\beta}}_j\delta^{6}}\\&
\le \frac{C}{\lambda^{\hat{\alpha}}_i\lambda^{3-\hat{\beta}}_j\delta^{\alpha}}\left(\frac{1}{\lambda_j \delta}\right)^{2\hat{\beta}-3}+\frac{C}{\lambda^{\hat{\alpha}}_i\lambda^{\hat{\beta}}_j\delta^{6}}.
\end{split}
\end{equation*}
Thus, we have for \;$I_2$
\begin{equation}\label{esti2int}
I_2\le \frac{C}{\lambda^{\hat{\alpha}}_i\lambda^{\hat{\beta}}_j\delta^{6}}.
\end{equation}
Next, for \;$I_1$\; we have 
\begin{equation*}
\begin{split}
I_1&=\int_{B_{a_i}(\delta)}\left(\frac{\lambda_i}{1+\lambda^2_i|x-a_i|^2}\right)^{\hat{\alpha}}\left(\frac{\lambda_j}{1+\lambda^2_jG^{-2}_{a_j}(x)}\right)^{\hat{\beta}}\\&
= \int_{B_{0}(\lambda_i\delta)} \left(\frac{1}{1+|y|^2}\right)^{\hat{\alpha}}\left[\frac{1}{\frac{\lambda_i}{\lambda_j}+\lambda_i\lambda_jG^{-2}_{a_j}\left(a_i+\frac{y}{\lambda_i}\right)}\right]^{\hat{\beta}}.
\end{split}\end{equation*}
Thus, if \;$\varepsilon^{-2}_{ij}\sim \frac{\lambda_i}{\lambda_j}$\; then 
\begin{equation*}
\begin{split}
I_1&\le C \varepsilon^{2\hat{\beta}}_{ij}\left[\left(\frac{1}{\lambda_i\delta}\right)^{2\hat{\alpha}-3}+C\right]\\&\le  C \varepsilon^{\beta}_{ij}.
\end{split}
\end{equation*}
If $\varepsilon^{-2}_{ij}\sim \lambda_i \lambda_j G^{-2}_{a_i}(a_j)$\; and \;$|a_i-a_j|\geq2\delta$, then we have 
\begin{equation*}
\begin{split}
I_1&\le C \left(\frac{1}{\lambda_i\lambda_j\delta^2}\right)^{\hat{\beta}}\left[\left(\frac{1}{\lambda_i\delta}\right)^{2\hat{\alpha}-3}+C\right]\\&
\le C \frac{1}{\delta^3}\left(\frac{1}{\lambda_i\lambda_j}\right)^{\hat{\beta}}\le C \frac{1}{\delta^3}\left[\left(\frac{1}{\lambda_i\lambda_jG^{-2}_{a_i}(a_j)}\right)^{\frac{1}{2}}\right]^{\beta}
\\&\le C \frac{1}{\delta^3}\varepsilon^{\beta}_{ij}.
\end{split}
\end{equation*}
Now, if $\varepsilon^{-2}_{ij}\sim \lambda_i \lambda_j G^{-2}_{a_i}(a_j)$ and $|a_i-a_j|<2\delta$ , then we get
\[I_1\le C\int_{B_{0}(\lambda_i\delta)} \left(\frac{1}{1+|y|^2}\right)^{\hat{\alpha}}\left[\frac{1}{\frac{\lambda_i}{\lambda_j}+\lambda_i\lambda_j\left|a_i+\frac{x}{\lambda_i}-a_j\right|^{2}}\right]^{\hat{\beta}}.\]
Next, we define 
\[B=\left\{\frac{1}{2}|a_i-a_j|\le \frac{|y|}{\lambda_i}\le2|a_i-a_j|\right\}\]
and have 
\[I_1\le C\int_{B}\left(\frac{1}{1+|y|^2}\right)^{\hat{\alpha}}\left[\frac{1}{\frac{\lambda_i}{\lambda_j}+\lambda_i\lambda_j\left|a_i+\frac{x}{\lambda_i}-a_j\right|^{2}}\right]^{\hat{\beta}}\]
\[\hspace{1.5cm}+C\int_{B_{0}(\lambda_i\delta)-B} \left(\frac{1}{1+|y|^2}\right)^{\hat{\alpha}}\left[\frac{1}{\frac{\lambda_i}{\lambda_j}+\lambda_i\lambda_j\left|a_i+\frac{x}{\lambda_i}-a_j\right|^{2}}\right]^{\hat{\beta}}.\]
For the second term, we have
\begin{equation*}
\begin{split}
\int_{B_{0}(\lambda_i\delta)-B} \left(\frac{1}{1+|y|^2}\right)^{\hat{\alpha}}\left[\frac{1}{\frac{\lambda_i}{\lambda_j}+\lambda_i\lambda_j\left|a_i+\frac{x}{\lambda_i}-a_j\right|^{2}}\right]^{\hat{\beta}}
&\le C \varepsilon^{\beta}_{ij}\left[\left(\frac{1}{\lambda_i\delta}\right)^{\alpha-3}+C\right]\\&\le C \varepsilon^{\beta}_{ij}.
\end{split}
\end{equation*}
For the first term, we have
\begin{equation*}
\begin{split}
&\int_{B}\left(\frac{1}{1+|y|^2}\right)^{\hat{\alpha}}\left[\frac{1}{\frac{\lambda_i}{\lambda_j}+\lambda_i\lambda_j\left|a_i+\frac{x}{\lambda_i}-a_j\right|^{2}}\right]^{\hat{\beta}}
\\&\le C\left(\frac{1}{1+\lambda^2_i|a_i-a_j|^2}\right)^{\hat{\alpha}}\int_{|y|\le 2\lambda_i|a_i-a_j|}\left[\frac{1}{\frac{\lambda_i}{\lambda_j}+\frac{\lambda_j}{\lambda_i}\left|y+\lambda_i(a_i-a_j)\right|^2}\right]^{\hat{\beta}}\\&
\le C\left(\frac{1}{1+\lambda^2_i|a_i-a_j|^2}\right)^{\hat{\alpha}}\int_{|z|\le 4\lambda_i|a_i-a_j|}\left[\frac{1}{\frac{\lambda_i}{\lambda_j}+\frac{\lambda_j}{\lambda_i}|z|^2}\right]^{\hat{\beta}}
\\&\le C\left(\frac{1}{\frac{\lambda_j}{\lambda_i}+\lambda_i\lambda_j|a_i-a_j|^2}\right)^{\frac{\alpha}{2}}\int_{|z|\le 4\lambda_j|a_i-a_j|}\left[\frac{1}{1+|z|^2}\right]^{\hat{\beta}}.
\end{split}
\end{equation*}
If $\lambda_j|a_i-a_j|$ is bounded, then we get
\begin{equation*}
\begin{split}
I_1&\le C\left(\frac{1}{\frac{\lambda_j}{\lambda_i}+\lambda_i\lambda_j|a_i-a_j|^2}\right)^{\frac{\alpha}{2}}\\&\le C\varepsilon^{\beta}_{ij}. 
\end{split}
\end{equation*}
If \;$\lambda_j|a_i-a_j|$\; is unbounded , then we get
\begin{equation*}
\begin{split}
I_1&\le C\left(\frac{1}{\frac{\lambda_j}{\lambda_i}+\lambda_i\lambda_j|a_i-a_j|^2}\right)^{\frac{\alpha}{2}}\left(\lambda_j|a_i-a_j|\right)^{3-2\hat{\beta}}\\
&\le C \left(\frac{1}{1+\lambda^2_i|a_i-a_j|^2}\right)^{\hat{\alpha}+\hat{\beta}-\frac{3}{2}}\left(\frac{\lambda_i}{\lambda_j}\right)^{\hat{\beta}}\\
&\le C \left(\frac{1}{\frac{\lambda_j}{\lambda_i}+\lambda_i\lambda_j|a_i-a_j|^2}\right)^{\hat{\beta}}\left(\frac{1}{1+\lambda^2_i|a_i-a_j|^2}\right)^{\hat{\alpha}-\frac{3}{2}}\\&\le C \varepsilon^{\beta}_{ij}.
\end{split}
\end{equation*}
Thus, we have for \;$I_1$
\begin{equation}\label{esti1f}
I_1\le \frac{C}{\delta^3}\varepsilon^{\beta}_{ij}.
\end{equation}
On the other hand, using the estimate for \;$I_2$ (see \eqref{esti2int}), we have
\begin{equation}\label{esti2f}
I_2=O\left(\frac{\varepsilon^{\beta}_{ij}}{\delta^6}\right).
\end{equation}
Hence, combining \eqref{esti1f} and \eqref{esti2f}, we have 
\[\int_\Omega u^{\alpha}_{a_i},_{\lambda_i} u^{\beta}_{a_j},_{\lambda_j}=O\left(\frac{\varepsilon^{\beta}_{ij}}{\delta^6}\right).\]
\end{pf}
\vspace{4pt}

\noindent
\section{Algebraic topological argument}
In this section, we present the algebraic topological argument for existence. We start by fixing some notation from algebraic topology.  
For  a topological space \;$Z$, $H_{*}(Z)$\; denotes the singular homology of \;$Z$ with \;$\Z_2$\; coefficients. If\;$Y$\; is a subspace of \;$Z$, then $H_*(Z, Y)$ stands for the relative homology with \;$\Z_2$\; coefficients of the topological pair \;$(Z, Y)$. For a map $f:Z\rightarrow Y$ \;with \;$Z$\; and \;$Y$\; topological spaces, \;$f_*$\; denotes the induced map in homology. If \;$f: (Z, Y)\longrightarrow (W, X)$ \; is a map with \;$(Z, Y)$\; and \;$(W, X)$\; topological pairs, then $ f_*$ denotes the induced map in relative homology. Furthermore, we discuss some algebraic topological tools needed for our application of the Barycentre technique of Bahri-Coron\cite{bc}. We start with the following observation. Since
 \;$\Omega$\; is a smooth domain of  \;$\R^3$\; which is non-contractble, then there exists \;$n\in \{1, 2\}$\; such that \;$H_n(\Omega)$\; is not trivial, see \cite{bc} (see page 1 just after Theorem 1). 
Hence, as  in \cite{bc} (see beginning of page 263), we have there exists \;$M$\; a smooth compact connected \; $n$-dimensional manifold without boundary and a continuous map \begin{equation}\label{conth}
h\; \;: \;M\longrightarrow \Omega
\end{equation}
 such that if we denote by \;$[M]$\; the class of orientation (modulo \;$2$\;) of \;$M$, then \;$h_*([M])\neq 0$. Moreover, we have clearly the existence of a compact smooth manifold with boundary \;$K_0$\; such that 
 \begin{equation}\label{k0}
 h(M)\subset K_0\subset \Omega.
 \end{equation}
We recall the space of formal barycenter of \;$M$\; that we need for our Barycenter technique for existence. For \;$p\in \N^*$,\; the set of formal barycenters of \;$M$\; of order \;$p$\; is defined as 
 \begin{equation}\label{eq:barytop}
B_{p}(M)=\{\sum_{i=1}^{p}\alpha_i\d_{a_i}\;:\;a_i\in M, \;\alpha_i\geq 0,\;\; i=1,\cdots, p,\;\,\sum_{i=1}^{p}\alpha_i=1\},\;
\;B_0(M)=\emptyset,
\end{equation}
where $\delta_{a}$\; for \;$a\in M$\; is the Dirac measure at \;$a$. we have the existence of  \;$\Z_2$\; orientation classes  (see \cite{bc})
\begin{equation}\label{orientation_classes}
w_p\in H_{np-1}(B_{p}(M), B_{p-1}(M)), \;\;\;\;\;p\in \N^*.
\end{equation}
Now to continue, we fix \;$\delta$\; small such that \;$0<2\d\leq \varrho_0$ \;where \;$\varrho_0$ \;as in \eqref{varrho} with \;$K$\; is replaced by \;$K_0$ \;and \;$K_0$\; is given by \eqref{k0}. Moreover, we choose \;$\theta_0>0$\; and smalll. After this, we let \;$\l$\; varies such that \;$0<\frac{1}{\l}\leq \theta_0\delta$\; and associate for every \;$p\in \N^*$\; the map
$$
f_p(\l): B_p(M)\longrightarrow H^{1, +}_0(\Omega) 
$$
defined by the formula
$$
f_p(\l)(\sigma)=\sum_{i=1}^p\alpha_i u_{h(a_i), \l}, \;\;\;\;\sigma=\sum_{i
=1}^p\alpha_i\d_{a_i},
$$
where\;$h$\; is as in \eqref{conth} and \;$u_{h(a_i), \l}$ is as \eqref{ual} (with \;$a$\; replaced by \;$h(a_i)$).
\vspace{6pt}

\noindent
As in Proposition 3.1in \cite{martndia2} and Proposition 6.3 in \cite{nss}, using Corollary \ref{interact2}, Corollary \ref{interact4}, Corollary \ref{sharpenergy}, Lemma \ref{interact5}, and  Lemma \ref{interact6}, we have the following multiple-bubble estimate.
\begin{pro}\label{eq:baryest}
There exist \;$\bar C_0>0$\; and \;$\bar c_0>0$\; such that for every \;$p\in \N^*$, $p\geq 2$ and every \;$0<\varepsilon\leq \varepsilon_0$, there exists \;$\l_p:=\l_p(\varepsilon)$ such that for every \;$\l\geq \l_p$ and for every $\sigma=\sum_{i=1}^p\alpha_i\delta_{a_i}\in B_p(M)$, we have
\begin{enumerate}
 \item
If \;$\sum_{i\neq j}\varepsilon_{i, j}> \varepsilon$\; or there exist \;$i_0\neq j_0$\; such that \;$\frac{\alpha_{i_0}}{\alpha_{j_0}}>\nu_0$, then
 $$
J_q(f_p(\l)(\sigma))\leq p^{\frac{2}{3}}\mathcal{S}.
$$ 
 \item
If \;$\sum_{i\neq j}\varepsilon_{i, j}\leq \varepsilon$\; and for every \;$i\neq j$\; we have \;$\frac{\alpha_{i}}{\alpha_j}\leq\nu_0$, then
$$
J_q(f_p(\l)(\sigma))\leq p^{\frac{2}{3}}\mathcal{S}\left(1+ \frac{\bar C_0}{\l}-\bar c_{0}\frac{(p-1)}{\l}\right),
$$
where \;$\varepsilon_{ij}$\; is as in \eqref{varepij} with \;$(a_i, a_j)$\; replaced by \;$(h(a_i), h(a_j))$ and $\l_i=\l_j=\l$, $\varepsilon_0$\; is as in \eqref{varepsilon0} and $\nu_0$ is as in \eqref{nu0}.
\end{enumerate}

\end{pro}
\vspace{6pt}

\noindent
As in Lemma 4.2  in \cite{martndia2} and Lemma 6.4 in \cite{nss}, we have the selection map \;$s_1$ (see \eqref{eq:mini}), Lemma \ref{deform} and Corollary \ref{sharpenergy} imply the following topological result.
\begin{lem}\label{eq:nontrivialf1}
Assuming that \;$J_q$\; has no critical points, then there exists \;$\bar \l_1>0$\;  such that for every \;$\l\geq \bar \l_1$, 
$$
f_1(\l)\; : \;(B_1(M),\; B_0(M))\longrightarrow (W_1, \;W_0)
$$
is well defined and satisfies
$$
(f_1(\l))_*(w_1)\neq 0\;\;\text{in}\;\;H_{n}(W_1, \;W_0).
$$
\end{lem}
\vspace{6pt}

\noindent
As in Lemma 4.3  in \cite{martndia2} and Lemma 6.5 in \cite{nss}, we have the selection map \;$s_p$ (see \eqref{eq:mini}), Lemma \ref{deform} and Proposition \ref{eq:baryest} imply the following recursive topological result.
\begin{lem}\label{eq:nontrivialrecursive}
Assuming that \;$J_q$\; has no critical points, then  there exists \;$\bar \l_p>0$\; such that for every \;$\l\geq \bar\l_p$, 
$$
f_{p+1}(\l): (B_{p+1}(M),\; B_{p}(M))\longrightarrow (W_{p+1}, \;W_{p})
$$
and 
$$
f_p(\l): (B_p(M), \;B_{p-1}(M))\longrightarrow (W_p, \; W_{p-1})
$$
are well defined and satisfy
$$(f_p(\l))_*(w_p)\neq 0\;\; \text{in}\;\; \;\;H_{np-1}(W_p, \;W_{p-1})$$ implies
$$(f_{p+1}(\l))_*(w_{p+1})\neq 0\;\; \text{in} \;\;H_{n(p+1)-1}(W_{p+1}, \;W_{p}).$$
\end{lem}
\vspace{6pt}

\noindent
Finally, as in Corollary 3.3 in \cite{martndia2} and Lemma 6.6 in \cite{nss}, we clearly have that Proposition \ref{eq:baryest}  implies the following result. 
\begin{lem}\label{eq:largep}
Setting \;$$\bar p_{0}:=[1+\frac{\bar C_{0}}{\bar c_{0}} ]+2$$ with \;$\bar C_0$\; and \;$\bar c_0$\; as in Proposition \ref{eq:baryest} and recalling \eqref{dfwp}, we have  there exists \;$\hat \l_{\bar p_0}>0$\; such that \;$\forall\l\geq\hat \l_{\bar p_0}$, 
$$
f_{\bar p_{0}}(\l)(B_{\bar p_{0}}(M))\subset W_{{\bar p}_0-1}.
$$
\end{lem}

\noindent
\begin{pfn} {of Theorem \ref{non-contractible}} \\
As in \cite{martndia2} and \cite{nss}, the theorem follows by a contradiction argument from Lemma \ref{eq:nontrivialf1} - Lemma \ref{eq:largep}.
\end{pfn}

\end{document}